\documentclass[onefignum,oneeqnum]{siamart171218}

\usepackage[utf8]{inputenc}
\usepackage[T1]{fontenc}
\usepackage{amsfonts}
\usepackage[mathscr]{euscript}
\usepackage{xfrac}
\usepackage{graphicx}
\usepackage{empheq}
\usepackage{caption}
\usepackage{subcaption}
\usepackage{url}
\usepackage{relsize}
\usepackage[normalem]{ulem}
\usepackage{algorithm}
\usepackage[noend]{algpseudocode}
\usepackage{enumitem}
\usepackage{placeins}
\usepackage{array}
\usepackage{float}
\usepackage{cellspace}
\colorlet{red}{black}
\colorlet{orange}{black}
\colorlet{purple}{black}

\usepackage{cite}
\newtheorem{remark}{Remark}
\crefname{subsection}{Section}{Sections}
\Crefname{subsection}{Section}{Sections}
\crefname{section}{Section}{Sections}
\Crefname{section}{Section}{Sections}
\crefrangeformat{figure}{Figures~#3#1#4 and~#5#2#6}
\Crefrangeformat{figure}{Figures~#3#1#4 through~#5#2#6}
\crefrangeformat{section}{Sections~#3#1#4 and~#5#2#6}
\Crefrangeformat{section}{Sections~#3#1#4 through~#5#2#6}
\crefrangeformat{subsection}{Sections~#3#1#4 and~#5#2#6}
\Crefrangeformat{subsection}{Sections~#3#1#4 through~#5#2#6}
\crefrangeformat{equation}{equations~(#3#1#4) through~(#5#2#6)}
\Crefrangeformat{equation}{Equations~(#3#1#4) through~(#5#2#6)}

\newcommand*\elide{\textup{[\,\dots]}}

\DeclareMathOperator{\sinc}{sinc}
\renewcommand{\i}{{\mathrm{i}}}
\renewcommand{\d}{{\mathrm{d}}}
\newcommand{\e}{{\mathrm{e}}}

\title{High-order, Dispersionless ``Fast-Hybrid'' Wave Equation
  Solver. Part~I: $\mathcal{O}(1)$ Sampling Cost via Incident-Field
  Windowing and Recentering}

\author{Thomas G. Anderson\thanks{Computing \& Mathematical Sciences, California Institute of
Technology (\email{tanderson@caltech.edu, obruno@caltech.edu}).}
\and Oscar P. Bruno\footnotemark[1] \and Mark Lyon\thanks{Department of Mathematics \& Statistics,
University of New Hampshire (\email{mark.lyon@unh.edu})}
}

\headers{``Fast-Hybrid'' Wave Equation Solver Pt. I.}{T. G. Anderson, O. P. Bruno, and M. Lyon}
\begin{document}

\date{}

\maketitle

\begin{abstract}
  This paper proposes a frequency/time hybrid integral-equation method
  for the time dependent wave equation in two and three-dimensional
  spatial domains. Relying on Fourier Transformation in time, the
  method utilizes a fixed (time-independent) number of
  frequency-domain integral-equation solutions to evaluate, with
  superalgebraically-small errors, time domain solutions for
  arbitrarily long times. The approach relies on two main elements,
  namely, 1)~A smooth time-windowing methodology that enables accurate
  band-limited representations for arbitrarily-long time signals, and
  2)~A novel Fourier transform approach which, in a time-parallel
  manner and without causing spurious periodicity effects, delivers
  numerically dispersionless spectrally-accurate solutions. A similar
  hybrid technique can be obtained on the basis of Laplace transforms
  instead of Fourier transforms, but we do not consider the
  Laplace-based method in the present contribution. The algorithm can
  handle dispersive media, it can tackle complex physical structures,
  it enables parallelization in time in a straightforward manner, and
  it allows for time leaping---that is, solution sampling at any given
  time $T$ at $\mathcal{O}(1)$-bounded sampling cost, for arbitrarily
  large values of $T$, and without requirement of evaluation of the
  solution at intermediate times. The proposed frequency-time
  hybridization strategy, which generalizes to any linear partial
  differential equation in the time domain for which frequency-domain
  solutions can be obtained (including e.g.\ the time-domain Maxwell
  equations), and which is applicable in a wide range of scientific
  and engineering contexts, provides significant advantages over other
  available alternatives such as volumetric
  discretization{\color{red}, time-domain integral equations,} and
  convolution-quadrature approaches.
\end{abstract}

\begin{AMS}
    65M80, 65T99, 65R20
\end{AMS}

%--------------------------------------------------------------------------------------------------
%--------------------------------------------------------------------------------------------------
\section{Introduction} \label{Introduction}
%--------------------------------------------------------------------------------------------------
%--------------------------------------------------------------------------------------------------

This paper, Part~I of a two-part contribution, proposes a fast
frequency-time hybrid integral-equation method for the solution of the
time domain wave equation in two- and three-dimensional spatial
domains. Relying on 1)~A smooth time-windowing methodology that
enables accurate band-limited representations for arbitrary long time
signals, and 2)~A novel FFT-accelerated Fourier transform strategy
(which, without requiring finer and finer meshes as time grows, is
amenable to time parallelism and does not give rise to spurious
periodicity effects), the proposed approach delivers numerically
dispersionless solutions with numerical errors that decay faster than
any power of the frequency mesh-size used. For definiteness, the
theoretical discussions in the present paper are restricted to
configurations for which the time-dependent excitations propagate
along a single incidence direction---which is, in fact, one of the
most common incident field arising in applications{\color{red}---but
  our numerical-results section includes examples that incorporate
  incident fields of other (generic) types
  (\Cref{tab:3D_banjai_comparison} and \Cref{fig:3D_wideband_sphere}).
  The development of algorithms for treatment of the general-incidence
  case on the basis of precomputation strategies that utilize plane
  waves or other bases of incident fields will be left for future work
  (but see \Cref{rem:pwe}). Similarly, general window tracking
  strategies based on characteristics of field-decay in two and
  three-dimensional configurations (including treatment of trapping
  structures), and parallel implementations exploiting the
  methodology's inherent parallelism in space and time, will be
  presented elsewhere.}

In practice the proposed methodology enjoys a number of attractive
properties, including high accuracy without numerical dispersion
error; an ability to effectively leverage existing frequency-domain
scattering solvers for arbitrary, potentially complex spatial domains;
an ability to treat dispersive media; dimensional reduction (if
integral equation methods are used as the frequency domain solver
components); natural parallel decoupling of the associated
frequency-domain components; and, most notably, time-leaping, time
parallelism, and $\mathcal{O}(1)$ cost for solution sampling at
arbitrarily-large times without requirement of intermediate time
evaluation.  A similar hybrid technique can be obtained on the basis
of Laplace transforms instead of Fourier transforms. Use of the
Laplace-based technique would be advantageous for treatment of certain
types of initial/boundary-value problems with non-vanishing initial
conditions, but we do not consider a Laplace-based approach in any
detail in the present contribution. We also note that the ideas
inherent in the proposed fast-hybrid approach may be applied, more
generally, to any time-domain problem whose frequency-domain
counterpart can be treated by means of an efficient frequency-domain
approach.

A wide literature exists, of course, for the treatment of the
classical wave equation problem.  Among the many approaches utilized
in this context we find finite-difference and finite-element time
domain methods~\cite{Taflove:00, Lee:97} (FDTD and FETD,
respectively), retarded potential boundary integral equation
methods~\cite{HaDuong:86,HaDuong:03, Yilmaz:04, Epstein:16},
Huygens-preserving treatments for odd-dimensional spatial
domains~\cite{Petropavlovsky:18}, and, most closely related to the
present work, two hybrid frequency-time methodologies, namely, the
Laplace-transform/finite-difference convolution quadrature
method~\cite{Lubich:94, Banjai:09, Banjai:10, Banjai:11, Banjai:12, Banjai:14,
  Betcke:17}, and the Fourier-transform/operator-expansion
method~\cite{Mecocci:00}.  A brief discussion of the
%method~\cite{Mecocci:00, Recchioni:03}.  A brief discussion of the
character of these methodologies is presented in what follows.

The FDTD approach and related finite-difference methods underlie most
of the wave-equation solvers used in practice. In these approaches the
solution on the entire spatial domain is obtained via finite
difference discretization of the PDE in both space and time. For the
ubiquitous exterior-domain problems, the use of absorbing boundary
conditions is necessary to render the problem computationally
feasible---which has in fact been an important and challenging problem
in itself~\cite{Hagstrom:09,Berenger:94,Engquist:77,Bayliss:80}. Most importantly, however, finite-difference methods
suffer from numerical dispersion, and they therefore require the use
of fine spatial meshes (and, thus, fine temporal meshes, for
stability) to produce accurate solutions.  Numerical dispersion errors
therefore present a significant obstacle for high frequency and/or
long time simulations via methods based on finite-difference spatial
discretizations. FETD methods provide an additional element of
geometric generality, but they require creation of high-quality finite
element meshes (which can be challenging for complex three-dimensional
structures). Further, like FDTD methods, they entail use of absorbing
boundary conditions, and they also generally give rise to detrimental
dispersion errors (also called ``pollution errors'' in this
context~\cite{Babuska:97}).

Integral-equation formulations based on direct discretization of the
time-domain retarded-potential Green's function, on the other hand, require treatment of the Dirac
delta function and thus give rise to integration domains given by the
intersection of the light cone with the overall scattering
surface~\cite{HaDuong:86,HaDuong:03}. These approaches generally
result in relatively complex overall schemes for which it has proven
rather challenging to ensure stability~\cite{Epstein:16}, and
which have typically been implemented in low-order accuracy setups
and, thus, with significant numerical dispersion error. Accelerated
versions of these methods have also been
proposed~\cite{Yilmaz:04}. {\color{purple} Motivated by the work
in~\cite{Epstein:16}, temporally and spatially
high-order time domain integral equation schemes have recently been
proposed~\cite{Barnett:19}.} Loosely related to this class of methods
are recent work on discretizations which are Huygens-preserving---that
is, treatments of the retarded potential operators with the advantage
that they do not entail an increasing amount of computational work for
increasing time, at least in odd dimensions~\cite{Petropavlovsky:18}.

Hybrid time-frequency methods rely on transform techniques to evaluate
time domain solutions by synthesis from sets of frequency domain
solutions; clearly the necessary solutions of (decoupled)
frequency-domain problems can be obtained via parallel
computation. The Convolution Quadrature (CQ) method~\cite{Lubich:94}
is a prominent example of this class of approaches. This method relies
on the combination of a finite-difference time discretization and a
Laplace transformation to effectively reduce the time domain wave
equation to a set of modified Helmholtz equations over a range of
frequencies.  There has additionally been some interest in the direct
use of Fourier transformations in
time{\color{red}~\cite{Mecocci:00,Douglas:93}} to decouple
the
%transformations in time~\cite{Mecocci:00,Recchioni:03} to decouple the
time-domain problem into frequency-domain sub-problems.  In detail,
assuming a Gaussian-modulated incident time-pulse, the
approach~\cite{Mecocci:00} evaluates Fourier integrals on
%approach~\cite{Mecocci:00,Recchioni:03} evaluates Fourier integrals on
the basis of a Gauss-Hermite quadrature rule, and it obtains the
necessary frequency-domain solutions by means of a certain ``operator
expansion method''{\color{red}; earlier efforts in a single spatial
dimension~\cite{Douglas:93} recognized the advantages of hybrid methods for simulating
wave phenomena in complex, and specifically attenuating
media, but relied on a low-order accurate midpoint quadrature rule for Fourier
integrals}.  In {\color{red} all cases}, however, the number of
frequency-domain solutions required, and hence the associated
computing and memory costs, grow linearly with the number $N$ of time
steps used to evolve the solution to a given final time $T$. (A more
detailed discussion of previous hybrid methods, including CQ and
direct-transform methods, is presented
in~\cref{sec:prev_hybrid_work}.)

To address these difficulties, for incident pulses of arbitrary
duration the proposed approach employs a \emph{smoothly time-windowed
  Fourier transformation technique} (detailed in
\cref{sec:time_partition}), which, without resorting to use of refined
frequency discretizations, re-centers the solution in time and thus
effectively handles the fast oscillations that occur in the scattering
solution as a function of the Fourier-transform variable
$\omega$. Therefore, in contrast to the CQ and the direct-transform
methods, the new approach can be applied in the presence of arbitrary
incident fields on the basis of a fixed set of frequency-domain
solutions. Other favorable properties of the method include its
time-parallel character, its time-leaping abilities and
$\mathcal{O}(1)$ cost of evaluation at any given time, however large,
and, therefore, its $\mathcal{O}(N)$ cost for a total full $N$
timestep history of the solution. The algorithm remains uniformly
(spectrally) accurate in time for arbitrarily long times, with
complete absence of temporal dispersion errors.

The proposed hybrid method relies on the use of a sequence of smooth
windowing functions (the sum of all of which equals unity) to smoothly
partition time into a sequence of windowed time-intervals.  The
claimed overall $\mathcal{O}(N)$ time cost with uniform-accuracy for
arbitrarily large times can be achieved for any given incident field
through use of time partitions of relatively large but fixed width,
leading to fixed computational cost per partition for arbitrarily long
times. In order to achieve such large-time uniform accuracy at fixed
cost per window, in turn, a new quadrature method for the evaluation
of windowed Fourier transform integrals is introduced which does not
require use of finer and finer discretizations for large
times---despite the increasingly oscillatory character, as time grows,
of a certain complex exponential factor in the transform integrands.
The time evaluation procedure requires computation of certain ``scaled
convolutions'' (with a $\sinc$ function kernel) which can be additionally
accelerated on the basis of the Fractional Fourier
Transform~\cite{BaileySwarztrauberSIAMRev:91}.

The hybrid methodology described in the present contribution lends
itself naturally, in a number of ways, to high-performance
load-balanced parallel computing. While full development of such
efficient parallelization strategies will be left for future work,
here we present some considerations in these regards. The simplest
parallel acceleration strategy in the context of the proposed method
concerns the set of frequency domain solutions it requires, which can
clearly be produced in an embarrassingly parallel fashion---whereby
frequency-domain problems are distributed among the available
computing cores. The evaluation of near fields, on the other hand,
also presents significant opportunities for parallel
acceleration. Indeed, the time-trace calculations on a prescribed
region $\mathcal{R}$ in space could be handled by distributing subsets
of $\mathcal{R}$ among various core groupings, or by relying on
frequency parallelization for evaluation of the necessary
frequency-domain near fields, or a combination of the two---depending
on 1)~The parallelization method used (if any) for the
frequency-domain problems
themselves~\cite{Bruno:01,BrunoGarza2018,Coifman:93,Bleszynski:96},
2)~The physical extent of the region $\mathcal{R}$, and 3)~The number
of frequencies that need to be considered for a given problem.  Time
parallelism, finally, can easily be achieved as a by-product of the
smooth time-partitioning approach. The multiple levels of parallelism
inherent in the algorithm should provide significant flexibility for
parallel implementations that exploit the differing capabilities of
various computer architectures.

This paper is organized as follows. After certain necessary
preliminaries are presented in \cref{sec:preliminaries} (including
well-known frequency domain integral formulations of the wave-equation
problem), the main components of the proposed approach are taken up in
\crefrange{sec:time_partition}{sec:fast_four_algs}. Thus
\cref{sec:time_partition} introduces the smooth time-partitioning
technique that underlies the proposed accelerated treatment of signals
of arbitrary long duration, while \cref{sec:fast_four_algs} puts forth
a new quadrature rule for the fast spectral evaluation of Fourier
transform integrals, with high-order accuracy and $\mathcal{O}(1)$
large-time sampling costs. An overall algorithmic description is
presented in \cref{sec:overall_algorithm}, and a variety of numerical
results, followed by some concluding remarks, finally, are presented
in \crefrange{sec:numerical_results}{sec:conclusion}. We believe that,
in view of its spectral time accuracy, absence of stability
constraints, fast algorithmic implementations, easy use in conjunction
with any existing frequency-domain solver, and bounded memory
requirements, the proposed method should prove attractive in a number
of contexts in science and engineering.

%--------------------------------------------------------------------------------------------------
%--------------------------------------------------------------------------------------------------
\section{Preliminaries}\label{sec:preliminaries}
%--------------------------------------------------------------------------------------------------
%--------------------------------------------------------------------------------------------------

%--------------------------------------------------------------------------------------------------
\subsection{{\color{orange}Wave equations in time domain and frequency domain}}\label{sec:diff_int_form}
%--------------------------------------------------------------------------------------------------
We consider the initial boundary value problem
\begin{subequations}\label{w_eq}
    \begin{align}
        \frac{\partial^2 u}{\partial t^2}(\mathbf{r}, t)& - c^2\Delta u
            (\mathbf{r}, t) = 0,\quad\mathbf{r} \in \Omega,\label{w_eq_a}\\
        u(\mathbf{r},0) &= \frac{\partial u}{\partial t}(\mathbf{r}, 0)
            = 0\\
        u(\mathbf{r}, t) &=  b(\mathbf{r}, t)
            \quad\mbox{for}\quad(\mathbf{r},t)\in\Gamma\times [0,T^\textit{inc}],\label{w_eq_c}
    \end{align}
\end{subequations}
for the time domain wave equation in the exterior domain
$\Omega\subset\mathbb{R}^d$ (the complement of a bounded set) for
$d=2,3$. The boundary of $\Omega$, which we will denote by $\Gamma$,
is an arbitrary Lipschitz surface for which an adequate frequency
domain integral-equation solver on $\Gamma$, or some alternative
frequency-domain method in the exterior of $\Gamma$, can be used to
solve the required frequency domain problems in the domain
$\Omega$. For definiteness, throughout this paper we assume a boundary
condition of the form~\cref{w_eq_c}, but similar treatments apply in
presence of boundary conditions of other types. Given an
incident field $u^\textit{inc}$, the selection $b = -u^\textit{inc}$
corresponds to a sound soft boundary condition for the total field
$u+u^\textit{inc}$ on the boundary of the scatterer:
\begin{equation*}
    u^\textit{tot}(\mathbf{r}, t) = u^\textit{inc}(\mathbf{r}, t) + u(\mathbf{r}, t)
        = 0,\quad\mathbf{r} \in \Gamma.
\end{equation*}
The Fourier transforms $U^t$ and $B^t$ of the solutions $u$ and $b$ of the wave
equation~\cref{w_eq} satisfy the Helmholtz problem with linear dispersion
relation $\kappa = \kappa(\omega) = \omega / c$,
\begin{subequations}\label{helmholtz}
    \begin{align}
        \Delta U^t(\mathbf{r}, \omega) + \kappa^2(\omega)U^t(\mathbf{r}, \omega) &= 0, \quad \mathbf{r} \in
        \Omega\\
        U^t(\mathbf{r}, \omega) &= B^t(\mathbf{r}, \omega), \quad \mathbf{r} \in \Gamma.
    \end{align}
\end{subequations}

\begin{remark}\label{rem:pwe}
  For definiteness, in the present contribution we restrict
  {\color{purple}most of our discussion} to one of the most common
  incident-field functions {\color{purple} $b = b(\mathbf{r}, t)$}
  arising in applications, namely, incident fields impinging along a
  single direction $\mathbf{p}$:
    \begin{equation}\label{single_incid_field}
        b(\mathbf{r}, t) = \frac{1}{2\pi} \int_{-\infty}^\infty B^t(\omega)
        \e^{\i \frac{\omega}{c} (\mathbf{p}\cdot\mathbf{r} - ct)}\,\d\omega{\color{purple}\quad\mbox{and}\quad B^t(\omega) = \int_{-\infty}^\infty a(t) \e^{\i \omega t}\,\d t,}
    \end{equation}
    for some {\color{purple} compactly supported function $a(t)$
      (consistent with the time interval of interest
      in~\cref{w_eq_c})}. {\color{purple} Note that in particular,
      $b(\mathbf{r}, t) = - u^{inc}(\mathbf{r}, t) = a(t - \mathbf{p}
      \cdot \mathbf{r}/c)$.}  {\color{red} In order to ensure the re-usability of
    the required set of frequency-domain solutions (see
    Section~\ref{sec:time_partitioning_hybrid}), arbitrary-incidence
    fields could either be treated by means of source- or
    scatterer-centered spherical expansions; or synthesis relying on
    principal-component analysis, etc. Such methodologies will be
    considered elsewhere, but a brief discussion in this regard is
    presented in \Cref{sec:pwe1}.  }
\end{remark}

\begin{remark}\label{rem:pwe2}
   The super-index $t$ in $B^t(\omega)$ (\Cref{single_incid_field}),
   indicates that the variable $\omega$ in this function's argument is
   the Fourier variable corresponding to $t$. In general, for any
   given function $f(\mathbf{r}, t)$ (resp. $f(t)$), $F^t(\mathbf{r},
   \omega)$ (resp. $F^t(\omega)$) will be used to denote the partial
   (resp.\ full) temporal Fourier transform of $f$ with respect to
   $t$, as indicated e.g.\ in \Cref{transform_pair}. Although only
   partial Fourier transforms in time are used in the present Part~I,
   the notation is adopted here to preserve consistency with
   Part~II---in which partial- or full-transforms with respect to both
   temporal and spatial variables are used.
\end{remark}

\subsection{Previous hybrid methods: Convolution
  quadrature~\cite{Lubich:94} and direct Fourier transform in
  time~\cite{Mecocci:00}}\label{sec:prev_hybrid_work}
%--------------------------------------------------------------------------------------------------
%--------------------------------------------------------------------------------------------------

As mentioned in~\cref{Introduction}, two hybrid time-domain methods
(i.e., methods that rely on transformation of the time variable by
means of Fourier or Laplace transforms) have previously been proposed,
namely, the Convolution Quadrature
method~\cite{Lubich:94,Banjai:09,Banjai:10,Banjai:11,Banjai:12,Banjai:14,Betcke:17}
and the direct Fourier transform
method~\cite{Mecocci:00}. The Convolution Quadrature
%method~\cite{Mecocci:00,Recchioni:03}. The Convolution Quadrature
method employs a discrete convolution that is obtained as temporal
finite-difference schemes are solved by transform methods.  Like the
method introduced in the present contribution, in turn, the direct
Fourier transform method is based on direct Fourier synthesis of
time-harmonic solutions. The following two sections briefly review
these two methodologies.

%--------------------------------------------------------------------------------------------------
\subsubsection{Convolution quadrature}\label{sec:prev_hybrid_CQ}
%--------------------------------------------------------------------------------------------------
The CQ algorithms result as the Z-transform is applied to the forward
recurrence relation arising from finite-difference temporal
semi-discretizations of the problem~\cref{w_eq}. A key point is that the
resulting time domain solution is \emph{itself an approximation of the
chosen temporal finite-difference approximation of the solution}{\color{red}.}
%Equivalently, the accuracy of the CQ method can at best
%be expected to match the accuracy of the underlying time
%discretization used.
In detail, utilizing the Z-transform, a finite-difference time
discretization of the wave equation can be reformulated as a set of
modified Helmholtz problems. The discrete time domain solution is then
obtained by evaluation of the inverse Z-transform of the frequency
domain solutions by means of trapezoidal-rule quadrature.
(References~\cite{Lubich:94,Betcke:17} provide further elaboration on
the connections of the CQ method to Z-transforms and convolutions,
respectively.)  As a result, the solutions produced by this method
accumulate temporal and spatial discretizations errors at each
timestep as well as overall inversion errors arising from the
approximate quadrature used in the inversion of the $Z$-transform. The
reliance of the CQ algorithm on a certain ``infinite-tail'' in the
time-history presents certain difficulties also. A brief discussion of
the character of these approximation methods is presented in what
follows.

The characteristics of a particular implementation of the CQ algorithm
is determined by the choice made for time-domain finite-difference
discretization, the spectral character of the discrete
frequency-domain solver used~\cite{Betcke:17}, and the methods
utilized for numerical inversion of the $Z$-transform. Existing CQ
approaches have primarily utilized the second-order accurate BDF2 time
discretization~\cite{Banjai:09}, but recent work~\cite{Banjai:14}
proposes the use of higher-order $m$-stage Runge-Kutta schemes. In all
cases the numbers of required frequency-domain solutions (which equals
$N_f$ for single-stage methods and $mN_f$ for $m$-stage methods, where
$N_f$ denotes the number of frequencies used to invert the $Z$-transform),
grows in a roughly linear fashion with the size of the
time interval for which the solution is to be produced. Thus, the cost
of the $m$-stage CQ approaches is $\mathcal{O}(mN_t)$, where $N_t$
denotes the number of time-steps taken. Stability and accuracy
considerations presented in~\cite{Betcke:17}, further, do suggest that
the stability of the CQ algorithm may be linked to certain
``scattering poles'' of the spatial solution operator which depend on
both the geometry of the spatial domain and the choice of the
frequency-domain formulation used. Reference~\cite{Betcke:17} further
suggests that the error of the contour integral discretization in the
CQ method (which is typically effected via the trapezoidal rule) can
dominate the error in the overall CQ time-stepping algorithm (even
under the well established $N_f = N_t$ setup), and that this
difficulty can be mitigated by overresolving the problem in
frequency domain---that is, using $N_f > N_t$. Furthermore, approximation
errors in FFT-accelerated evaluation of the Cauchy
integral formula for required weights in the $Z$-transform inversion typically
imply~\cite[\S 3.3]{Banjai:12} a maximum achievable overall accuracy of
$\sqrt{\varepsilon_\mathrm{mach}} \approx 10^{-8}$, where
$\varepsilon_\mathrm{mach}$ denotes double precision machine epsilon. Finally, reduction of
order of \textit{temporal} convergence is observed at points in the
near-field as an observation point approaches the scatterer~\cite{Banjai:11}.

In addition to $Z$-transform-inversion errors, the numerical
dissipation and dispersion introduced by the underlying time-domain
finite difference discretizations present an additional important
source of error in the CQ approach~\cite{Chen:12,Banjai:12}. These
errors can be managed by utilizing a number of timesteps which varies
superlinearly with frequency~\cite[\S 4.3]{Banjai:12} (that is, faster
than the number of sampling points required for uniformly accurate
interpolation), but the computational cost associated with such
procedures can be significant.

The memory requirements of the CQ method can be significantly impacted
by its reliance on a certain ``infinite time-tail'', which is
described e.g.\ in~\cite[Chapter~5]{Sayas:16}:
\begin{quote}
  The sequence of problems~\elide~presents the serious disadvantage of
  having an infinite tail.  In other words, the passage through the
  Laplace domain introduces a regularization of the wave equation that
  eliminates the Huygens’ principle that so clearly appears in the
  time domain retarded operators and potentials.
\end{quote}
The infinite tail impacts on the computing costs of the CQ method in
two different ways, namely: 1)~As the CQ timestep tends to zero for a
fixed final time $T$; and 2)~As the final time $T$ grows for a fixed
timestep. While the growth in point~1) can be slowed to certain extent
by appealing to Laplace-domain decay rates of compactly-supported,
smooth incident data~\cite{Banjai:09} (whose $Z$-transform counterpart
generally decays much faster than the error arising from the
time-stepping scheme utilized, and can thus be broadly neglected up to
the prescribed error tolerance), the infinite-tail growth in point~2)
has remained untreated, and it does give rise to linear growth in the
overall CQ computing and memory cost per timestep as $T\to\infty$.

\subsubsection{Direct Fourier transform in time}\label{sec:prev_hybrid_direct_fourier}
Without reliance on finite difference approximations, the direct
Fourier transform method proposed in~\cite{Mecocci:00} proceeds by
Fourier transformation of the time domain wave equation followed by
solution of the resulting Helmholtz equations for a range of
frequencies and inverse transformation to the time-domain using the
transform pair
\begin{equation}\label{transform_pair}
    U^t(\mathbf{r}, \omega) =\int_{-\infty}^{\infty} u(\mathbf{r}, t) \e^{\i
    \omega t}\, \d t,\qquad
    u(\mathbf{r}, t) = \frac{1}{2\pi} \int_{-\infty}^\infty U^t(\mathbf{r},
    \omega) \e^{-\i \omega t}\, \d \omega
\end{equation}
(see \Cref{rem:pwe2}). In detail, for general boundary values $b =
b(\mathbf{r}, t)$ (\Cref{w_eq}), reference~\cite{Mecocci:00} uses a
plane wave representation of the form
\begin{align}\label{plane_wave_fourier_pair}
    B^t(\mathbf{r}, \omega) &= \frac{1}{(2\pi)^d}  \int_{S^{d-1}}
    B^{\mathbf{r}, t}(\mathbf{p}, \omega) \e^{\i \kappa(\omega) \mathbf{p} \cdot \mathbf{r}}
    \,\d \mathbf{p},
\end{align}
so that full solution $U^t$ can be reconstructed on the basis of the
solutions $U^t = U^t_\mathbf{p}$ of Helmholtz
problems~\cref{helmholtz}, where
$\kappa = \kappa(\omega) = \omega / c$, and where the sound-soft
boundary values are given by the plane wave
$\e^{\i\kappa(\omega)\mathbf{p}\cdot\mathbf{r}}$ in the direction of
the vector $\mathbf{p}$.  The numerical examples in~\cite{Mecocci:00}
assume overall boundary {\color{orange}data} of this form for a single incidence
vector $\mathbf{p}$---that is, a uni-directional incident wave field.

Importantly, the resulting direct Fourier method does not suffer from
dispersion errors in the time variable. In the
contribution~\cite{Mecocci:00} the needed Helmholtz solutions are
obtained by means of a certain ``operator-expansion'' technique, and
assumes the incident field is given by a plane wave modulated by a
Gaussian envelope in frequency domain in order that the needed Fourier integrals
are approximated using the classical Gauss-Hermite quadrature rule.

Except for simple geometries, the use of the operator-expansion method
limits the overall accuracy to the point that in many cases it is
difficult to discern convergence. This difficulty could be addressed
by switching to a modern, more effective, frequency-domain
technique. Most significantly, however, the use of any generic
numerical integration procedure for the evaluation of the necessary
inverse Fourier transforms for large $t$, including the highly
accurate Gauss-Hermite rule used in~\cite{Mecocci:00}, does lead to
difficulties---in view of the highly-oscillatory character, with
respect to $\omega$, of the exponential factor in the right-hand
expression in~\cref{transform_pair} for large values of $t$. Indeed,
evaluation of the aforementioned inverse Fourier transform for required
time sample values up to the final time $T$ on the basis
of such procedures requires use of a number $N$ of frequency
discretization points (and hence a number of required frequency-domain
solutions) which is proportional to $T$. Calling $P$ the average cost of these
frequency-domain solutions, and including the overall
$\mathcal{O}(T^2)$ computational cost required by the evaluation of
the $\mathcal{O}(T)$-cost Gauss-Hermite inverse transform for each
$t_n\leq T$, the overall cost of the algorithm~\cite{Mecocci:00} can
be estimated as
\begin{equation}
TP + T^2.
\end{equation}
This estimate must be contrasted with the cost required by
classical finite-difference methods---which is proportional to the
first power of $T$.

\subsubsection{Accuracy and computational costs}
Estimates on the accuracy of hybrid methods follow from
well-established results on convergence of the associated
frequency-domain solution techniques together with corresponding
accuracy estimates on the underlying treatment of frequency/time
discretizations. The CQ method (\Cref{sec:prev_hybrid_CQ}) has
typically used low-order Galerkin spatial
discretizations~\cite{Banjai:10}, although the recent
contribution~\cite{Labarca:19} does incorporate a high-order
frequency-domain solver. Frequency/time discretization errors in the
CQ method, on the other hand, arise from the time-stepping scheme used
and the numerical discretization selected of a certain complex contour
integral. Typically, BDF2 is chosen as the underlying CQ time-stepping
scheme, yielding second order accuracy in time---but see
also~\cite{Banjai:14} for use of higher-order temporal CQ
discretizations and their associated computing costs. Concerning the
CQ complex-contour quadrature, on the other hand, the trapezoidal
quadrature rule that is used most often in this context can be an
important source of numerical error (see~\cite{Betcke:17} and the
previous discussion in~\cref{sec:prev_hybrid_CQ}). The direct Fourier
Transform method~\cite{Mecocci:00}
(\Cref{sec:prev_hybrid_direct_fourier}), in turn, exhibits high-order
Gauss-Hermite convergence in time, but generally poor spatial
convergence for the frequency-domain problems (but see
\Cref{sec:prev_hybrid_direct_fourier} in these regards). The fast
hybrid method proposed in the present paper, finally, relies on
well-known Nystr\"om frequency-domain methods, which generally exhibit
superalgebraically fast convergence (that is, convergence faster than
any power of the discretization mesh), together with exponentially
convergent methods for evaluating frequency/time transforms (except,
in the 2D cases with low-frequency content for which arbitrarily high
but not exponential convergence is obtained). In sum, the accuracy of
the CQ methods is mostly limited by the errors arising in the
time-stepping evolution scheme if the necessary complex integrations
are performed with sufficient accuracy. The direct Fourier
method~\cite{Mecocci:00} and the fast hybrid method proposed in this
paper, in turn, enjoy highly favorable convergence properties as
discretizations are refined.

The total computational costs required by the various hybrid
algorithms under consideration will be quantified in terms of the
number $N$ of time-points $t_n$ ($1\leq n \leq N$) at which the
solution is desired, as well as the average computing cost $P$
required by each one of the necessary frequency-domain solutions.
Roughly speaking (up to logarithmic factors), the CQ methods entail a
computing cost proportional to $NP$---that is, the method requires a number of
frequency-domain solutions that grows linearly with time. In more
detail, for example, reference~\cite[Sec. 4]{Banjai:09} proposes a CQ
algorithm for which it reports a computing cost of
$\mathcal{O}(N\log^2(N) P)$ operations.  According to
\Cref{sec:prev_hybrid_direct_fourier}, in turn, the Direct Fourier
Transform method requires $\mathcal{O}(N^2)+\mathcal{O}(N P)$
operations. In contrast, as shown in \Cref{sec:fast_four_algs}, the
proposed fast hybrid method requires $rP + \mathcal{O}(N)$ operations
to evaluate the solution at $N$ time points (where $r$ is the number,
independent of $N$, of frequency-domain solutions required by the
method to reach a given accuracy for arbitrarily long time).

Note that while the previous hybrid methods require the solution of an
increasing number of Helmholtz problems as time grows, the proposed
fast hybrid method does not---a fact which lies at the heart of the
method's claimed $\mathcal{O}(1)$-in-time sampling cost for
arbitrarily large times $t$.  In terms of memory storage, the fast
hybrid method requires $r V$ memory units for sampling at arbitrarily
large times, where $V$ denotes the average value of the storage needed
for each one of the necessary frequency-domain solutions. Of course,
storage of the entire time history on a given set of spatial points,
which may or may not be desired, does require a total of
$\mathcal{O}(N)$ memory units.

%--------------------------------------------------------------------------------------------------
\subsection{Frequency-domain representation}\label{sec:cfie_repr}
%--------------------------------------------------------------------------------------------------
The method of layer potentials provides an effective technique for the
solution of the frequency-domain problem~\cref{helmholtz}. The
layer-potential method we use in this paper relies on the
frequency-domain single and adjoint double-layer operators,
\begin{equation*}
    (S_\omega\psi)(\mathbf{r}) = \int_\Gamma G_{\omega}(\mathbf{r}, \mathbf{r'})
    \psi(\mathbf{r'})\,\d\sigma(\mathbf{r'}), \quad\mbox{and}\quad
    (K^*_\omega\psi)(\mathbf{r}) = \int_\Gamma \frac{\partial G_{\omega}(\mathbf{r},
    \mathbf{r'})}{\partial n(\mathbf{r})} \psi(\mathbf{r'})\,\d\sigma(\mathbf{r'}),
\end{equation*}
for $\mathbf{r} \in \Gamma$, {\color{orange} where $G_\omega$ denotes the
fundamental solution of the Helmholtz equation at frequency $\omega/c$
---which in the two- and three-dimensional cases is given,
respectively, by
\begin{equation}
    G_\omega(\mathbf{r}, \mathbf{r}') = \frac{\i}{4}
    H_0^{(1)}(\frac{\omega}{c}|\mathbf{r} - \mathbf{r}'|) \quad\mbox{and}\quad
    G_\omega(\mathbf{r}, \mathbf{r}') = \frac{\e^{\i\frac{\omega}{c}
    |\mathbf{r} - \mathbf{r}'|}}{4\pi|\mathbf{r} - \mathbf{r}'|}.
\end{equation}
Use of the {\color{purple}Green's} function $G_\omega$ and Green's third
identity yields the frequency domain field representation
    \begin{equation}\label{freq_domain_repr}
      U^t(\mathbf{r}, \omega) = \int_\Gamma \psi^t(\mathbf{r'}, \omega) G_{\omega}
        (\mathbf{r}, \mathbf{r}')\,\d\sigma(\mathbf{r}'),\quad \mathbf{r} \in
        \mathcal{R} \subset \Omega.
    \end{equation}
    where $\psi^t$, which equals the boundary values of the normal
    derivative of the total field ($\psi^t(\mathbf{r},\omega) =
    \frac{\partial U^{t,\textit{tot}}(\mathbf{r},
      \omega)}{\partial n(\mathbf{r})}$){\color{purple},} may be obtained} as the
    solution of the direct integral equation
\begin{equation}\label{single_layer_int_eq}
    (S_\omega \psi^t)(\mathbf{r}, \omega) = B^t(\mathbf{r}, \omega), \quad \mathbf{r} \in \Gamma.
\end{equation}
Further discussion of the connections between hybrid frequency/time
formulations and time-domain integral representations will be
presented in Part~II.

Unfortunately, equation~\cref{single_layer_int_eq} is not uniquely
solvable for certain values of $\omega$.  Making use of the auxiliary
adjoint double-layer operator we obtain the \emph{uniquely solvable} direct
combined field integral equation formulation
({\color{red}see~e.g.~\cite{ChandlerWilde:12}}):
\begin{equation}\label{CFIE_direct}
    \frac{1}{2}\psi^t(\mathbf{r}, \omega) {\color{red}+} (K^*_\omega \psi^t)(\mathbf{r}, \omega) -
    \i\eta(S_\omega \psi^t)(\mathbf{r}, \omega)= \frac{\partial B^t(\mathbf{r}, \omega)}{\partial
    n(\mathbf{r})} {\color{red}-} \i\eta B^t(\mathbf{r}, \omega), \quad \mathbf{r} \in \Gamma.
\end{equation}
A wide literature exists for the numerical solution of boundary
integral equations of this type. In this paper we use Nystr\"om
methods to discretize and solve the integral
equations~\cref{CFIE_direct} for all desired frequencies. In the case
$d = 2$ (resp.\ $d =3$) the Nystr\"om method described
in~\cite[\S 3.5]{ColtonKress} (resp.\ in~\cite{BrunoGarza2018}) is used.

%-------------------------------------------------------------------------------
\section{Smooth time-partitioning Fourier-transformation
  strategy}\label{sec:time_partition}
%-------------------------------------------------------------------------------

An efficient smooth time-partitioning {\em
  ``windowing-and-recentering''}\/ solution algorithm is proposed in
this section which is based on a number of novel methodologies. The
algorithm first expresses the solution $u$ of~\cref{w_eq}, for
arbitrary large times $T$, in terms of solutions $u_k$ arising from
incident fields that are compactly supported in time: $u(\mathbf{r},
t) = \sum_{k=1}^K u_k(\mathbf{r}, t)$ ($K = \mathcal{O}(T)$). Assuming
the incident fields can be represented with a given error tolerance
$\varepsilon$ within a time-frequency bandwidth $W$, the re-centering
component of the strategy presented in
\Cref{sec:time_partitioning_hybrid} produces all of the
functions $u_k$ in terms of a certain \emph{fixed finite set}
$\mathcal{F} = \{\psi_\mathbf{p}^t(\cdot, \omega_j),\, (1 \leq j \leq J)\}$ of
frequency-domain solutions appropriate for the assumed temporal
bandwidth $W$ (cf.\ equations \cref{ukW} and \cref{Ukslow_repr}).
The re-utilization of a fixed set of boundary integral densities
$\left\{\psi^t_\mathbf{p}\right\}$ (and hence the requirement of a
fixed number $J$ of solutions of the integral
equation~\cref{CFIE_direct} for evaluation of $u(\mathbf{r},t)$ for
arbitrarily long times $t$), is a key element leading to the
effectiveness of the proposed algorithm for incident signals of
arbitrarily-long duration.

%-------------------------------------------------------------------------------
\subsection{Time partitioning, windowing and re-centering, and the
  Fourier Transform}
\label{sec:time_partition_gen}
%-------------------------------------------------------------------------------
Motivated by \Cref{transform_pair}, let $(f,F)$ denote a
Fourier Transform pair
\begin{equation}
  F(\omega) =\int_{0}^{T} f(t) \e^{\i \omega t}\,\d t, \quad f(t)  =
    \frac{1}{2\pi} \int_{-\infty}^\infty F(\omega) \e^{-\i \omega t}\,\d\omega,
\label{transform_pair_bounded}
\end{equation}
for a (finitely or infinitely) {\em smooth compactly supported}\/
function $f(t)$, assumed zero except for $t \in [0,T]$ ($T>0$) (as
there arise, e.g., in the smooth time-partitioning strategy described
in \Cref{sec:time_partitioning_hybrid}). In this case the
Fourier transform on the left-hand side
of~\cref{transform_pair_bounded} is an integral over a finite (but
potentially large) time interval.

\begin{figure}
    \centering
    \includegraphics[width=0.49\textwidth]{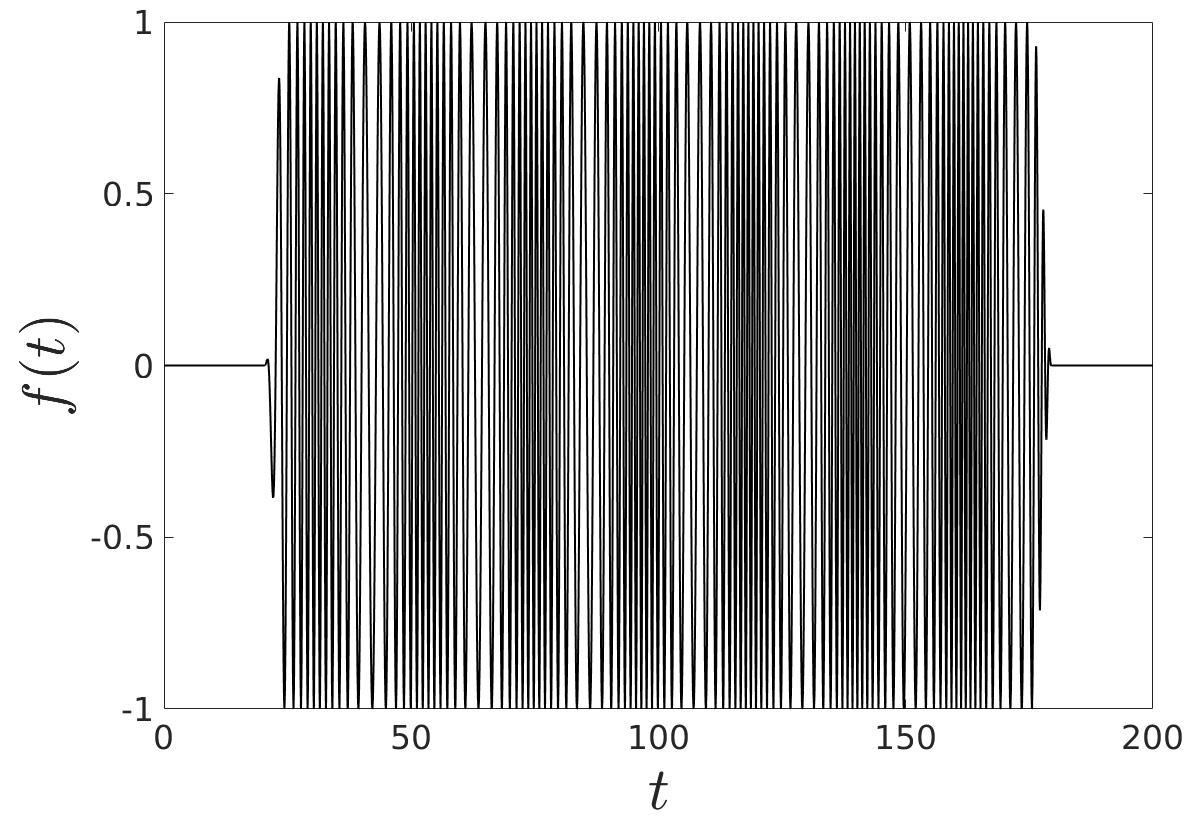}
    \includegraphics[width=0.49\textwidth]{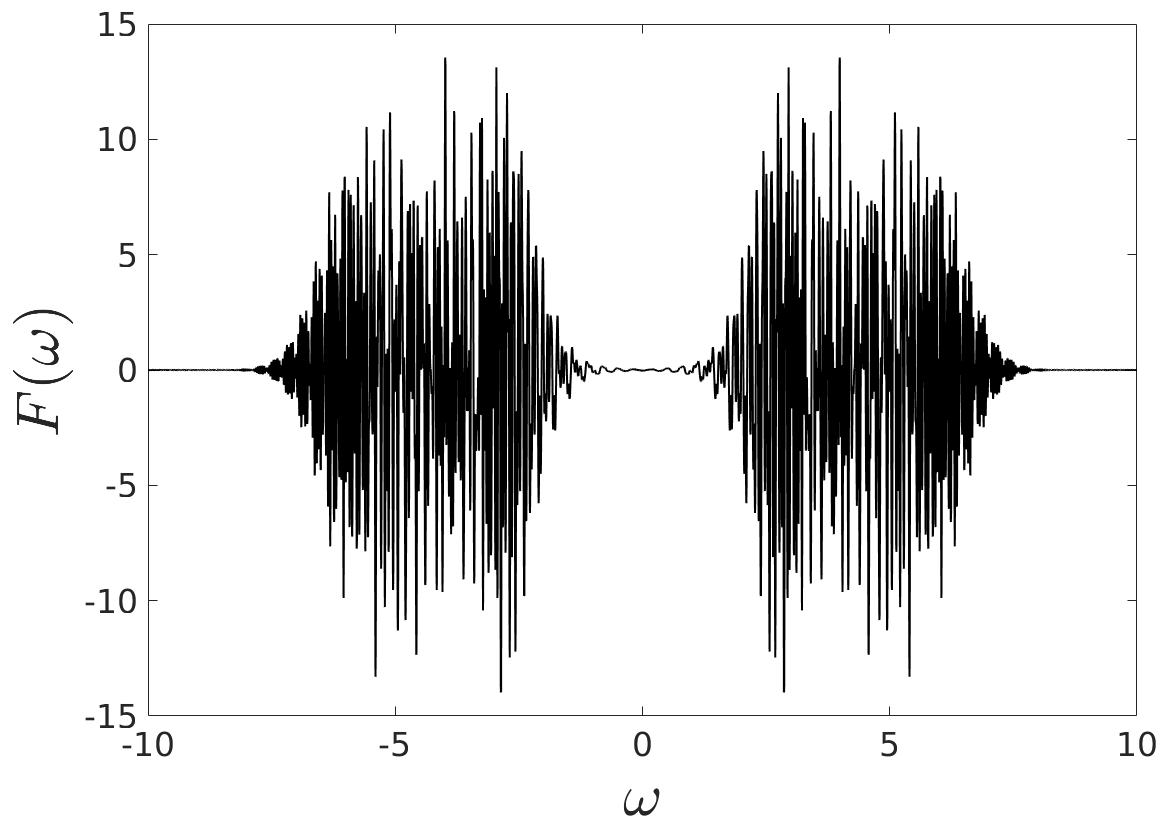}
    \caption{Left: Smooth, long duration time signal $f(t)$ {\color{red} as
    given in \cref{chirp_incidence}, windowed to have support}
      in the interval $20\leq t \leq 180$. Right: Real part of the
      Fourier Transform $F(\omega)$ of $f(t)$. The Fourier transform
      $F(\omega)$ is highly oscillatory on account of the large $t$
      values contained in the left-hand integrand
      in~\cref{transform_pair_bounded}.}
    \label{fig:Ms1}
\end{figure}

In the context of our problem it is useful to consider the dependence
of the oscillation rate of the function $F(\omega)$ on the parameter
$T$.  \cref{fig:Ms1} demonstrates the situation for a representative
``large-$T$'' chirped function $f$ depicted on the left-hand image in
the figure: the Fourier transform $F(\omega)$, depicted on the
right-hand image is clearly highly oscillatory. Loosely speaking, the
highly-oscillatory character of the function $F(\omega)$ stems from
corresponding fast oscillation in the factor $\e^{\i\omega t}$
contained in the left-hand integrand in \Cref{transform_pair_bounded}
for each fixed large value of $t$.  The consequence is that a very
fine discretization mesh $\omega_j$, containing $\mathcal{O}(T)$
elements, would be required to obtain $f(t)$ from $F(\omega)$ on the
basis of the right-hand expression
in~\cref{transform_pair_bounded}. In the context of a hybrid
frequency-time solver, this would entail use of a number
$\mathcal{O}(T)$ of applications of the most expensive part of the
overall algorithm: the boundary integral equations solver---which
would make the overall time-domain algorithm unacceptably slow for
long-time simulations. This section describes a new Fourier transform
algorithm that produces $f(t)$ (left image in \cref{fig:Ms1}) within a
prescribed accuracy tolerance, and for any value of $T$, however
large, by means of a $T$-independent (small) set of discrete frequency
values $\omega_j$ ($-W\leq \omega_j\leq W$, $j=0,\dots, J$).

The proposed strategy for the large-$T$ Fourier transform problem is
based on use of a partition-of-unity (POU) set $\mathcal{P} =
\{w_k(t)\,|\,k=1,\dots,K\}$ of ``well-spaced'' windowing functions,
where $w_k$ is supported in a neighborhood of the point $s=s_k$ for
certain ``support centers'' $s_k\in [0,T]$ ($1\leq k\leq K$)
satisfying, for some constants $H_1,H_2>0$, the minimum-spacing
property $s_{k+1} - s_k \geq H_1$, as well as the maximum width
condition $w_k(t) = 0$ for $|t-s_k|> H_2$ and the partition-of-unity
relation $\sum_{k=1}^K w_k =1$. Setting $H = H_1 = H_2$ in our test
cases we use POU sets based on the following parameter selections
\begin{enumerate}[label=\alph*)]
    \item $s_{k+1} - s_k =  3H/2$,
    \item $w_k(t)=1$ in a neighborhood $|t - s_k| < H/2$,
    \item $w_k(t)=0$ for $|t- s_k|>H$, and
    \item\label{three} $\sum_{k=1}^K w_k(t) = 1\quad \mbox{for all}\quad t\in
        [0,T]$.
\end{enumerate}
Note that, since $H$ is (or, more generally $H_1$ and $H_2$ are)
$T$-independent, the integer $K$ is necessarily an $\mathcal{O}(T)$
quantity. In practice we use the prescription $w_k(t) = w(t - s_k)$, with
partition centers and window function given by $s_k = 3(k-1)H/2$ {\color{red}
(the parameter choice $H = 10$ was used in all cases in this article)} and
\begin{equation}\label{window_func}
    {\color{purple}w(t)} =
    \begin{cases}
        1 - \eta({\color{red}\frac{t + H}{H/2}}), & -H \le t \le -H/2\\
        1, & -H/2 < t < H/2\\
        \eta(\frac{t - H/2}{H/2}), & H/2 \le t \le H\\
        0, & |t| > H,
    \end{cases}
\end{equation}
respectively, where we use the smooth windowing function $\eta \in
C_c^\infty([-1, 1])$, $\eta(u) = \exp(\frac{2e^{-1/u}}{u-1})$.

%{\color{red},
%though we emphasize that this is by no means the only window function proposed
%for Fourier Transforms: see~\cite{Harris:78} for windowing in a discrete
%setting and also the related literature of the Short Time Fourier
%Transform~\cite{Grochenig:01}}.
\begin{figure}
    \centering

    \includegraphics[width=0.49\textwidth]{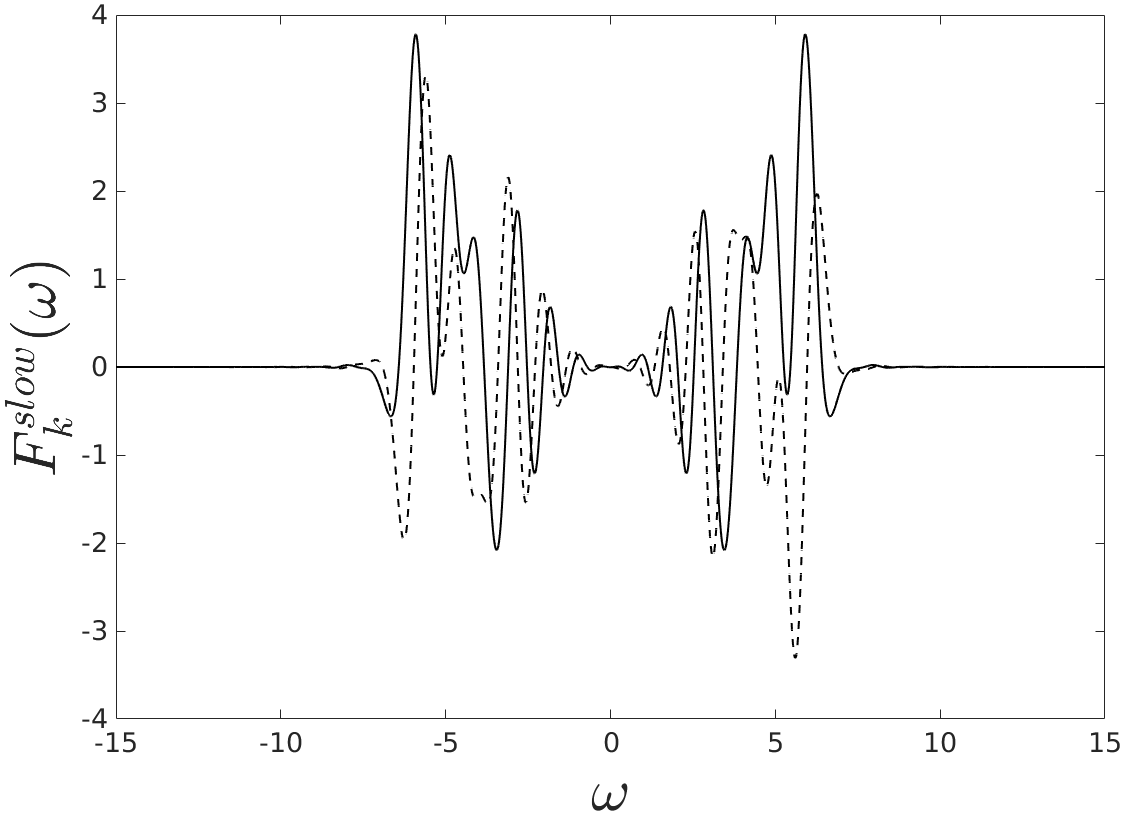}
    \includegraphics[width=0.49\textwidth]{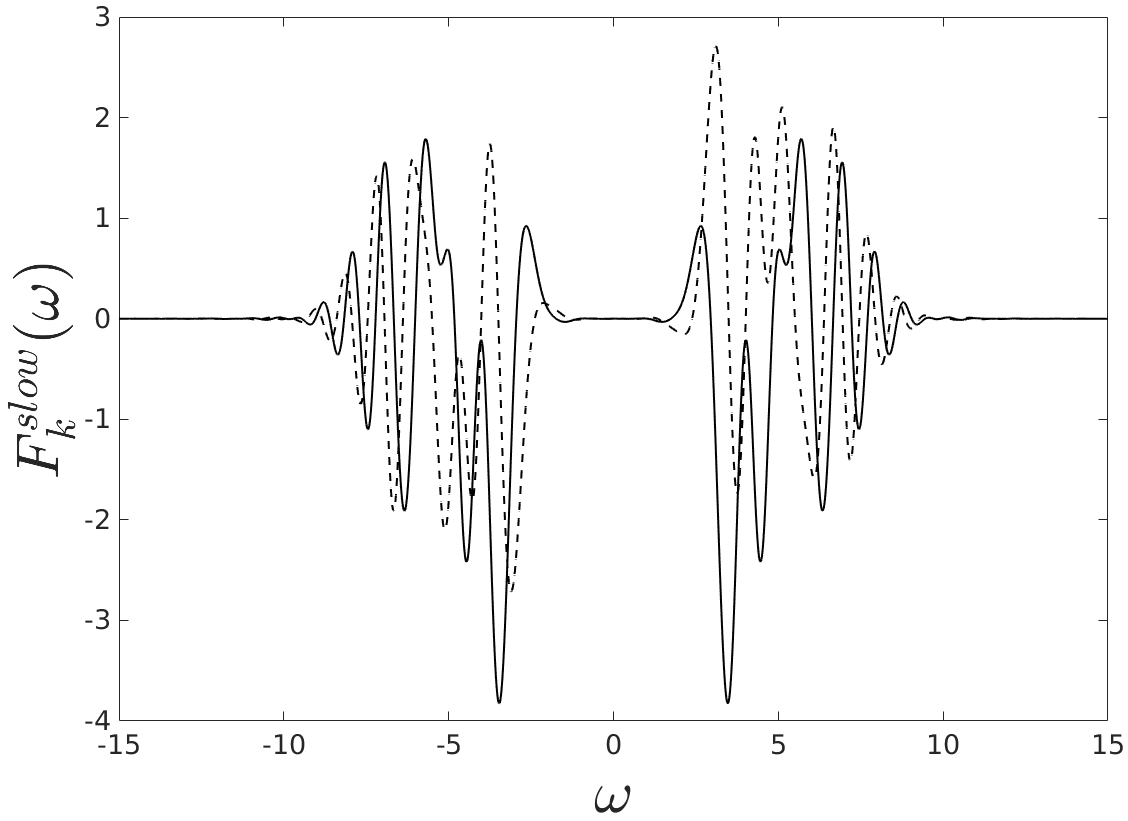}
    \caption{Fourier Transform of two windowed partitions of the long
      duration signal shown in \Cref{fig:Ms1}, each with partition
      width $H=10$. With reference to the text, the left and right
      figures depict the transform corresponding, respectively, to
      partition centers at $s_k = 35$ ($k = 4$) and $s_k = 155$
      ($k = 16$). In each case, the solid and dashed traces depict the
      real and imaginary parts of the Fourier Transform,
      respectively. The transforms are both less than $10^{-4}$
      outside the plotted region.}
    \label{fig:Ms2}
\end{figure}

Using the partition of unity $\mathcal{P}$ and letting $f_k(t) =
w_k(t) f(t)$, for $\omega\in[-W,W]$ we obtain the expression
\begin{equation}\label{Fk_sum}
  F(\omega) = \sum_{k=1}^K F_k(\omega),\quad\mbox{where}\quad
  F_k(\omega)=\int_{s_k-H_2}^{s_k+H_2} f_k(t) \e^{\i \omega t}\,\d t,
\end{equation}
{\color{red} which resembles the type of integrals used in connection with the windowed Fourier transform~\cite{Grochenig:01}.}
Now, centering the integration interval around the origin we obtain
\begin{equation}\label{fk_def}
F_k(\omega)= \int_{-H_2}^{H_2} f_k(t+s_k) \e^{\i \omega (t+s_k)}\,\d t=  \e^{\i \omega s_k} F_k^\textit{slow}(\omega)
\end{equation}
where
\begin{equation}\label{Fkslow_int}
  F_k^\textit{slow}(\omega)= \int_{-H_2}^{H_2} f_k(t+s_k) \e^{\i \omega t}\,\d t.
\end{equation}
The ``slow'' superscript refers to the fact that, since $t$
in~\cref{Fkslow_int} is ``small'' (it satisfies $-H_2\leq t \leq
H_2$), it follows that the integrand~\cref{Fkslow_int} only contains
slowly oscillating exponential functions of $\omega$, and thus
$F_k^\textit{slow}(\omega)$ is itself slowly
oscillatory. Thus~\cref{fk_def} expresses $F_k(\omega)$ as product of
two terms: the (generically) highly oscillatory exponential term $\e^{i \omega s_k}$
(which arises for signals whose support center is away from the origin \emph{in time}), on one hand, and the
slowly oscillatory term $F_k^\textit{slow}(\omega)$, on the
other. \cref{fig:Ms2} displays the real and imaginary parts of
$F_k^\textit{slow}$ for two values of $k$, namely $k=4$ and
$k=16$. Note that despite the differing centers in time, the functions
are similarly oscillatory and both are much less oscillatory than the
Fourier transform depicted in \cref{fig:Ms1}.

\begin{remark}\label{rem:1}
  Since $f$ is smooth and compactly supported, iterated integration by parts in
  the integral expressions that define $F(\omega)$, $F_k(\omega)$ and
  $F_k^\textit{slow}(\omega)$ (equations~\cref{transform_pair_bounded},
  \cref{fk_def} and~\cref{Fkslow_int}), and associated expressions for the
  derivatives of these functions of any positive order, shows that these
  functions and their derivatives decay {\color{red} as $1/\omega^n$ as
  $\omega\to\pm \infty$ for all $n > 1$ for which $f \in C^n$}. In other
  words, {\color{red} for smooth functions $f$} these three functions, along
  with each one of their derivatives respect to $\omega$ (of any order),
  decay superalgebraically fast as $\omega\to\pm \infty$. Additionally, in
  the two latter cases the superalgebraically fast decay (for each fixed
  order of differentiation) is uniform in $k$.
\end{remark}

\begin{remark}\label{rem:1a}
  Let $G:\mathbb{R}\to \mathbb{C}$, $G=G(\omega)$, denote a function
  that decays super-algebraically fast, along with each one of its
  derivatives, as $\omega\to\infty$. Then, repeated use of integration
  by parts on the inverse Fourier transform expression $ g (t)=
  \frac{1}{2\pi} \int_{-\infty}^\infty G(\omega)\e^{-\i\omega
    t}\,\d\omega$ shows that the error in the approximation
\[
g (t) \approx \frac{1}{2\pi} \int_{-W}^W G(\omega)\e^{-\i\omega t}\,\d\omega
\]
decays super-algebraically fast as $W\to\infty$.
\end{remark}

\begin{remark}\label{rem:1b}
  Let $G = G(\omega)$ denote a superalgebraically-decaying function,
  as in \Cref{rem:1a}. Then, repeated use of integration by parts in
  the integral expressions for the Fourier coefficients $g_n$ shows
  that the expansion of $G$ as a $2W$-periodic Fourier series
\[
    G(\omega) \approx \sum_{n=-\infty}^\infty g_n \e^{\i {\color{red}n}\pi
  \omega/W},\quad -W\leq \omega\leq W,
\]
together with all of its derivatives, converge to $G(\omega)$ and its
respective derivatives uniformly and super-algebraically fast, as
$W\to\infty$, throughout the interval $[-W,W]$.
\end{remark}

%-------------------------------------------------------------------------------
\subsection{Windowed and re-centered wave equation and  solutions with slow $\omega$ dependence}
\label{sec:time_partitioning_hybrid}
%-------------------------------------------------------------------------------

In order to evaluate numerically the solution of the
problem~\cref{w_eq} we apply the smooth time partitioning strategy
developed in \Cref{sec:time_partition_gen} to the boundary-condition
function $b(\mathbf{r},t)$ in~\cref{w_eq_c} (as a function of $t$ for
each fixed value of $\mathbf{r}$). Thus, using the window functions
$w_k(t)$ described in the previous section, we define {\color{purple}
  two different types of windowed boundary-condition functions
  $b_k=b_k(\mathbf{r},t) $, namely 1)~The function
\begin{equation}\label{bkdef2}
  b_k(\mathbf{r}, t) = w_k(t)b(\mathbf{r}, t),
\end{equation}
which can be used for general incident fields $b(\mathbf{r}, t)$, as
well as, 2)~Specifically for incident fields of the
form~\cref{single_incid_field}, the function
\begin{equation}\label{bkdef}
  b_{k}(\mathbf{r}, t) =  \frac{1}{2\pi} \int_{-\infty}^\infty
  B^t_{k}(\omega) \e^{\i (\kappa(\omega) \mathbf{p} \cdot \mathbf{r} - \omega t)}
  \,\d\omega\quad\mbox{with}\quad \kappa(\omega) = \omega / c,
\end{equation}
where, letting 
\begin{equation}\label{ak_signal_def}
a_k(t) = w_k(t)a(t)
\end{equation}
we have set
\begin{equation}\label{bk_signal_def}
  B_k^t(\omega) = \int_{-\infty}^{\infty}a_k(t)\e^{\i\omega t}\,\d t = \int_{-\infty}^{\infty}w_k(t)a(t)\e^{\i\omega t}\,\d t.
\end{equation}
Note that both definitions of $b_k(\mathbf{r}, t)$ imply
$\sum_{k=1}^K b_k(\mathbf{r}, t) = b(\mathbf{r}, t)$.  In view of
\Cref{rem:pwe} the present article uses the boundary condition
function~\cref{bkdef}.}

Letting $u_k(\mathbf{r}, t)$ ($1\leq k\leq K$) denote the solution
to~\cref{w_eq} with boundary-condition function $b(\mathbf{r}, t)$
substituted by $b_k(\mathbf{r}, t)$, we clearly have
\begin{equation}\label{uk_def}
u(\mathbf{r}, t) = \sum_{k=1}^K u_k(\mathbf{r}, t).
\end{equation}
This expression is the basis of the time-domain solver proposed in
this paper.
\begin{remark}\label{rem:tracking}
  As discussed extensively in Part~II, in view of Huygens' principle
  in three dimensions, and a certain windowing reallocation strategy
  in two dimensions, a fixed, geometry-dependent, number $M$
  (independent of $K = \mathcal{O}(T)$) of solutions $u_k$ need to be
  included for any space-time evaluation region, irrespective of the
  time duration $T$ for which the solution is evaluated. The
  geometry-dependence of the parameter $M$ relates closely to the
  trapping character~\cite{Morawetz:75,Morawetz:77} of the underlying
  scattering geometry{\color{red}; while $M$ may grow large for slowly-decaying
    scattered fields, the finite-duration character of the incident field is
    guaranteed by this theory to yield an overall-bounded number of active
    partitions}. This is the basis of tracking strategies for
  identifying the ``active'' time-partition solutions, which will be
  presented elsewhere.  \Cref{fig:CKK_plane_wave_offangle_tracking},
  which displays the computed functions $u_k$ for $k=1,2,3$, and in
  which the solution for each partition is only plotted if it exceeds
  a certain tolerance anywhere in the entire domain of interest,
  illustrates, in a rudimentary fashion, some of the principles
  inherent in those strategies.
\end{remark}

Accurate numerical approximations of the solutions $u_k$ can be
produced as indicated in what follows. {\color{purple}Considering the boundary
condition function $b_k(\mathbf{r}, t)$ in \Cref{bkdef},} Equations~\cref{Fk_sum} and~\cref{Fkslow_int} yield
\begin{equation}\label{Bk_def}
  B^t(\omega) = \sum_{k=1}^K B^t_k(\omega),
    \quad\mbox{and}\quad B^t_k(\omega) = \e^{\i\omega s_k}
    B_k^\textit{slow}(\omega).
\end{equation}
Thus, denoting by $U^t_k(\mathbf{r}, \omega)$ and
$U_k^\textit{slow}(\mathbf{r}, \omega)$ the frequency domain solutions
of the problem~\cref{helmholtz} with $B^t$ replaced by
$B^t_k(\omega)\e^{\i \kappa(\omega) \mathbf{p}\cdot\mathbf{r}}$ and
$B_k^\textit{slow}(\omega)\e^{\i \kappa(\omega) \mathbf{p}\cdot\mathbf{r}}$,
respectively, we obtain the representations
\begin{equation}\label{uk_rep}
  u_k(\mathbf{r}, t) = \frac{1}{2\pi} \int_{-\infty}^\infty U^t_k(\mathbf{r},
    \omega)\e^{-\i\omega
    t}\,\d\omega = \frac{1}{2\pi} \int_{-\infty}^\infty U_k^\textit{slow}(\mathbf{r},
  \omega)\e^{-\i\omega(t - s_k)}\,\d\omega.
\end{equation}
Since $U_k^\textit{slow}$ is approximately band-limited (because
$B_k^\textit{slow}$ is, see \cref{rem:1}), it follows
from~\cref{uk_rep} and \cref{rem:1a} that $u_k(\mathbf{r}, t)$ can be
approximated by the strictly band-limited function $u_k^W$:
\begin{equation}\label{ukW}
    u_k(\mathbf{r}, t) \approx u_k^W(\mathbf{r}, t) = \frac{1}{2\pi} \int_{-W}^W U_k^\textit{slow}(\mathbf{r}, \omega)
    \e^{-\i\omega(t - s_k)}\,\d\omega
\end{equation}
with superalgebraically small errors (uniform in $\mathbf{r}$, $t$
and $k$) as the bandwidth $W$ grows.

\Cref{sec:fast_four_algs} presents a quadrature algorithm that, on the
basis of a finite set of frequencies $\mathcal{F} = \{\omega_j: \,
j=1,\dots,J\}$, approximates, with errors uniform-in-$t$ and decaying rapidly
as $J$ increases, the highly-oscillatory integral~\cref{ukW},
and thus produces the numerical approximation $u_k^{W,J} \approx u_k^W$, by means of spectral interpolation of the
slowly-varying quantity $U_k^\textit{slow}(\cdot, \omega)$ {\em with
respect to $\omega$}. The quantity $U_k^\textit{slow}$, in turn, is
dependent on frequency-domain incident data $B_k^\textit{slow}$ obtained (via
use of the numerical transform techniques introduced in \Cref{sec:fourQuad}
{\color{purple} for the function $a_k = w_k a$})
from the relation
%\begin{equation}\label{Bkslow_eqn}
%    B_k^\textit{slow}(\omega) = \e^{-\i \kappa(\omega) \mathbf{p} \cdot \mathbf{r}}\int_{-H}^H b_{k}(\mathbf{r}, t + s_k)
%    \e^{\i\omega t}\,\d t,
%\end{equation}
{\color{purple}
\begin{equation}\label{Bkslow_eqn}
    B_k^\textit{slow}(\omega) = \int_{-H}^H a_{k}(t + s_k)
    \e^{\i\omega t}\,\d t,
\end{equation}
}
and boundary integral ``scattering'' densities $\psi_\mathbf{p}^t$ {\color{purple}, where $\psi_\mathbf{p}^t$ are solutions of \Cref{CFIE_direct} with $B^t$ replaced by $\e^{\i \kappa(\omega) \mathbf{p} \cdot \mathbf{r}}$}.
Specifically, {\color{purple} defining}
\begin{equation}\label{psikslow_repr}
   \psi_k^\textit{slow}(\mathbf{r}', \omega) = B_k^\textit{slow}(\omega)
     \psi_\mathbf{p}^t(\mathbf{r}', \omega),
\end{equation}
for the time-partition-specific boundary density we have
\begin{equation}\label{Ukslow_repr}
\begin{split}
  U_k^\textit{slow}(\mathbf{r}, \omega) = \int_\Gamma
  \psi_k^\textit{slow}(\mathbf{r}', \omega)
  & G_{\omega}(\mathbf{r}, \mathbf{r}')\,\d\sigma(\mathbf{r}') \\
  & =B_k^\textit{slow}(\omega) \int_\Gamma
  \psi_\mathbf{p}^t(\mathbf{r}', \omega) G_{\omega}(\mathbf{r},
  \mathbf{r}')\,\d\sigma(\mathbf{r}').
\end{split}
\end{equation}
It follows that for a given bandwidth $W$ all the needed function
values $U_k^\textit{slow}(\mathbf{r}, \omega_j)$ can be produced in
terms of the fixed ($k$-independent, $W$-dependent) finite set $\Psi =
\{\psi^t_\mathbf{p}(\cdot, \omega_j): \, j=1,\dots,J\}$ of boundary
integral densities. The re-utilization of the fixed ($k$-independent)
set $\Psi$ of ``expensive'' integral densities is a crucial element
leading to the efficiency of the overall hybrid algorithm.

%-------------------------------------------------------------------------------
\section{FFT-based $\mathcal{O}(1)$-cost Fourier transform at large
times}\label{sec:fast_four_algs}
%-------------------------------------------------------------------------------
This section presents an effective algorithm for the numerical
evaluation of truncated Fourier integrals of the form
\begin{equation}
    F(\omega) = \int_{-H}^H f(t) \e^{\i\omega t}\,\d t\quad\mbox{and}\quad
    f(t) = \frac{1}{2\pi} \int_{-W}^{W}
    F(\omega) \e^{-\i \omega t}\,\d\omega
\label{four_trans_limited}
\end{equation}
(cf. Equations~\cref{Fkslow_int} and~\cref{ukW}), at arbitrarily large
evaluation arguments $t$ and $\omega$. Here, it is assumed that $f$ is
a smooth function of time $t\in \mathbb{R}$ which vanishes outside the
interval $[-H,H]$. Similarly, with the possible exception of a
inverse-logarithmic singularity of $F(\omega)$ at $\omega = 0$ for
certain two-dimensional applications (see \cref{sec:sing_quad}), $F$
is an infinitely smooth function for all frequencies $\omega$---which
is additionally superalgebraically small outside the interval
$[-W,W]$.  The case in which a singularity exists in
\Cref{four_trans_limited} at $\omega = 0$ is handled in
\Cref{sec:sing_quad} by utilizing a decomposition of the form
\begin{equation}\label{fourier_sing_decomp_discussion}
    f(t)  = \left( \int_{-W}^{-w_c} + \int_{-w_c}^{w_c} + \int_{w_c}^W
    \right) F(\omega) \e^{-\i\omega t}\,\d\omega,
\end{equation}
together with a specialized quadrature rule for the middle integral;
the function $F$ is smooth (though not necessarily periodic) in the
integration intervals $[-W, -\omega_c]$ and $[\omega_c, W]$.

Use of trapezoidal rule integration might appear
advantageous in these contexts, since, for such boundary-vanishing
integrands the trapezoidal rule exhibits superalgebraically fast
convergence (at least in the smooth $F$ case), and, importantly,
unlike the Gauss-Hermite rule used in~\cite{Mecocci:00}, it can be
efficiently evaluated by means of FFTs. However, as the evaluation
arguments $t$ or $\omega$ grow, the integrands
in~\cref{four_trans_limited} become more and more oscillatory.  Both
the Gauss-Hermite and the trapezoidal rule (and, indeed, any
quadrature rule based on standard interpolation techniques) require
use of finer and finer meshes to avoid completely inaccurate
approximations as the evaluation argument increases
(see~\cref{sec:time_partition_gen}). Failure to resolve this
difficulty would lead to a fundamental breakdown in the algorithm---as
it would be necessary for the scheme to produce an increasing number
of (expensive) boundary integral equation solutions, leading to
rapidly increasing costs, as evaluation times grow.

\begin{remark}\label{inv_fwd_cov}
  For definiteness, the presentation in this section is restricted to
  the right-hand integral in~\cref{four_trans_limited}; the
  corresponding algorithm for the left-hand integral is entirely
  analogous.
\end{remark}

\begin{remark}
  A direct examination of the trapezoidal approximation
\begin{equation}
    f(t) =
    \frac{1}{2\pi} \int_{-W}^{W}
    F(\omega) \e^{-\i t \omega}\,\d\omega \approx \frac{W}{2\pi m}\sum_{k=0}^{m-1} F(\omega_k)
    \e^{-\i t\omega_k} \quad \left(\omega_k = -W + k\, \Delta \omega\right)
\label{trap_rule_four_trans}
\end{equation}
shows that, as is well known, quadrature errors in the trapezoidal
quadrature rule for ``large'' $t$ manifest themselves as ``aliasing'',
that is, spurious periodicity in the $t$
variable~\cite{BaileySwarztrauberSIAMJSciComp:94,Frazer:84,Griffith:90,Lin:92}.
\end{remark}

The method proposed in the present \cref{sec:fast_four_algs} resolves
the difficulties mentioned in the last paragraph: it eliminates
aliasing errors without recourse to frequency mesh refinement, and it
evaluates (on the basis of FFTs) the time-domain solution in constant
computing time per temporal evaluation point---so that, as in
finite-difference time-marching algorithms, the overall cost per
timestep of the time propagation algorithm does not grow as time
increases.

%--------------------------------------------------------------------------------------------------
\subsection{Smooth $F(\omega)$: FFT-based reduction to ``scaled
convolution''}\label{sec:fourQuad}
%--------------------------------------------------------------------------------------------------

This section considers the problem of evaluation of Fourier integrals
similar to those in~\cref{four_trans_limited} (or, in 2D contexts, the
integrals with smooth integrands in
\Cref{fourier_sing_decomp_discussion}) under the assumption that the
functions $f$ and $F$ are infinitely smooth in the domain of
integration. In the context of \Cref{ukW} in dimension $d=3$ the
smoothness assumption on $F$ is always satisfied, as it is for
dimension $d=2$ provided that, e.g., $F(\omega) =
U_k^\textit{slow}(\mathbf{r}, \omega)\e^{\i\omega s_k}$ vanishes in a
neighborhood of $\omega=0$. The singular $d=2$ case is tackled
in~\cref{sec:sing_quad}.

The proposed smooth-$F$ approach proceeds by trigonometric-series
expansion of the integrand function $F$ followed by use of certain
``scaled convolutions'' introduced in \cref{sec:qfsx}; a fast
FFT-based algorithm for evaluation of such convolution-like quantities
is then described in \cref{sec:disc_conv_fast}.
%--------------------------------------------------------------------------------------------------
\subsubsection{Transform approximation via Fourier Series expansion\label{sec:qfsx}}
%--------------------------------------------------------------------------------------------------
In this Section we develop a quadrature rule for the general transform
integral
\begin{equation}\label{four_trans_general}
  I_a^b\,[F](t) = \int_a^b F(\omega) \e^{-\i\omega t}\,\d\omega,
\end{equation}
or, equivalently,
\begin{equation}\label{four_trans_centered}
  I_a^b\,[F] (t) =
  \e^{-\i \delta t} \int_{-A}^A F(\delta + \omega) \e^{-\i \omega t} \,\d\omega,\quad \mbox{where}\quad A = \frac{b-a}{2}\quad\mbox{and} \quad\delta = \frac{b+a}{2}.
\end{equation}
Although $F(\delta + \omega)$ may not be a periodic function of
$\omega$ in the integration interval $[-A,A]$, for a prescribed
positive even integer $M$ we utilize a trigonometric polynomial of the
form
\begin{equation}\label{fourier_approx}
  F(\delta + \omega) \approx \sum_{m=-\sfrac{M}{2}}^{\sfrac{M}{2}-1} c_m \e^{\i\frac{2\pi}{P}
    m\omega}
\end{equation}
of a certain periodicity $P$, that closely approximates $F(\delta +
\omega)$ for $\omega\in[-A,A]$.
\begin{remark}\label{Fourier-FC}
  As indicated below, in the context of this paper $F(\delta +
  \omega)$ is most often a smoothly periodic function in $[-A,A]$
  (with $A$ equal to the bandlimit $W$); in such cases we take $P=2A$
  and~\eqref{fourier_approx} is obtained as a regular Discrete Fourier
  Transform (DFT) in $[-A,A]$.  Exceptions do arise in certain
  two-dimensional situations (\cref{sec:sing_quad}) where $F(\delta +
  \omega)$ is smooth but not periodic in $[-A,A]$ (cf.\ the first and
  last integrals in \cref{fourier_sing_decomp_discussion}); in such
  cases an accurate Fourier approximation of a certain period $P\ne
  2A$ is obtained in our algorithm on the basis of the FC(Gram)
  Fourier Continuation method~\cite{Bruno:10,Amlani:16}.  In the
  periodic case the errors inherent in the
  approximation~\cref{fourier_approx} tend to zero super-algebraically
  fast (faster than any negative power of $M$~\cite[Lemma
  7.3.3]{Atkinson:09}, cf.\ also \Cref{rem:1b}), while the errors
  arising from the Fourier Continuation method used in the
  non-periodic case decay as a user-prescribed negative power of $M$.
\end{remark}

Substituting~\cref{fourier_approx} into~\cref{four_trans_centered} and
integrating term-wise yields the approximation
\begin{equation}\label{approx}
    \begin{split}
        I_a^b\,[F] (t) &\approx
            \e^{-\i\delta t} \sum_{m=-\sfrac{M}{2}}^{\sfrac{M}{2}-1} c_m \int_{-A}^A
                \e^{-\i\frac{2\pi}{P}(\alpha t - m)\omega}\,\d\omega \\
            &= \e^{-\i\delta t} \sum_{m=-\sfrac{M}{2}}^{\sfrac{M}{2}-1} c_m \frac{P}
                {\pi(\alpha t - m)}\sin\left(\pi\frac{2A}{P}(\alpha t - m)\right),
    \end{split}
\end{equation}
where we have set $\alpha = \frac{P}{2\pi}$. In view of~\cref{approx},
for a given {\em user-prescribed (!)}\/ equi-spaced time-evaluation grid
$\left\{t_n = n\Delta t\right\}_{n = N_1}^{N_2}$ we may write, letting
$\beta = \alpha\Delta t $,
\begin{equation}\label{discConv}
  I_a^b\,[F] (t_n ) \approx \e^{-\i\delta t_n} \sum_{m=-\sfrac{M}{2}}^{\sfrac{M}{2} - 1}
  c_m b_{\beta n - m}, \quad \mathrm{where}\quad  b_q \coloneqq
  2A\sinc\left( \frac{2A}{P} q \right).
\end{equation}
Note that, paralleling the fast Fourier series convergence in the
periodic case, equations~\cref{approx} and~\cref{discConv} provide
super-algebraically close approximations of $I_a^b[F](t)$. In the
non-periodic case, these equations provide a user-prescribed algebraic
order of accuracy. In either case, the errors in \cref{approx}
and~\cref{discConv} are uniform in $t$ and $n$, respectively: for a
given error tolerance $\varepsilon$ there exists an integer $M_0$
(independent of $t$ and $t_n$) such that, for all $M\geq M_0$, the
approximation errors in~\cref{approx} and~\cref{discConv} are less
than $\varepsilon$ for all $t\in\mathbb{R}$ and all relevant values
$t_n$, respectively.

\begin{remark}\label{fourier_points_smooth}
  It is useful to note that the aforementioned $t$- and
  $t_n$-independent errors in \cref{approx} and~\cref{discConv} stem
  solely from corresponding errors in the expansion
  \cref{fourier_approx}---and thus, can be achieved on the basis of
  values of the function $F(\delta + \omega)$ on a fixed
  ($t$-independent) finite set $\mathcal{F}^\textrm{smooth}$ of
  frequency mesh points, cf.  \Cref{sec:overall_algorithm}.
\end{remark}

Since generically $\beta \neq 1$ {\color{red} (indeed, $\beta \not\in
\mathbb{Z}$ generically)}, the quantity $\sum_m c_m b_{\beta n
  - m}$ in~\cref{discConv} is not a discrete convolution, but it is,
rather, a ``discrete scaled convolution''~\cite{Nascov:09}.  Like
regular discrete convolutions, scaled convolutions can accurately be
produced by means of FFTs~\cite{Nascov:09}---although the algorithm
for scaled convolutions is somewhat more complicated than the standard
FFT convolution approach. Still, the fast scaled convolution algorithm is a
useful tool: it runs in $\mathcal{O}(L\log L)$ operations (where $L =
\mathrm{max}(N_2 - N_1,M)$) and it produces highly accurate results;
details are presented in~\cref{sec:disc_conv_fast}.

%-------------------------------------------------------------------------------
\subsubsection{FFT accelerated evaluation of scaled discrete convolutions}\label{sec:disc_conv_fast}
%-------------------------------------------------------------------------------
The quadrature method introduced in~\cref{sec:qfsx} reduces the
evaluation of the right-hand transform in~\cref{four_trans_limited}
for values $t = t_n$ (for a given range $0\leq n - n_0\leq N-1$ with
$n_0\in\mathbb{Z}$ and $N\in\mathbb{N}$) to evaluation of scaled
convolutions of the form
\begin{equation}\label{conv_sum}
  d_n = \sum_{m=-\sfrac{M}{2}}^{\sfrac{M}{2}-1} c_m b_{\beta m - \gamma
    n},\quad 0\leq n - n_0\leq N-1,
\end{equation}
where the coefficients $c_m$ are complex numbers that make up a
certain ``input vector'' $\vec{c} = (c_{-M/2},\dots,c_{-M/2-1})$, and
where the ``convolution kernel'' $b$ is a function of its real-valued
sub-index $q$: $b_q=b(q)$. (Compare~\cref{discConv}
and~\cref{conv_sum} and note the specific scaled convolution kernel
$b_q$ and parameter value $\gamma=1$ used in the former equation.)
This section presents an algorithm which evaluates the
sum~\cref{conv_sum} for all required values of $n$ at FFT speeds.

To describe the algorithm, let $L$ denote a certain positive even integer, to
be defined below, which is larger than or equal to the maximum of $N$ and $M$.
The convolution input vector $\vec{c}$ is symmetrically zero-padded to form a
new vector $\vec{c}=(c_{-L/2},c_{-L/2+1},\dots,c_{L/2-1})$ of length $L$. New
elements are also added to the list of evaluation indices in~\cref{conv_sum} so
that the overall list contains the $L$ elements in the indicial vector $\vec{n}
=(n_0-L/2,n_0-L/2+1,\dots, n_0 + L/2 - 1)$. Following~\cite{Nascov:09}, for
technical reasons the length $L$ is determined by the relation $L\geq L_0$,
where $L_0$ denotes the smallest even integer for which the kernel index
parameter $q = \beta m - \gamma n$ lies in the range $-L_0/2 \le q \le L_0/2 -
1$ for $-M/2 \le m \le M/2 - 1$, $0 \le n - n_0 \le N - 1$.  (As pointed out
below, selections satisfying $L>L_0$ are occasionally necessary to achieve a
prescribed error tolerance.) In view of these selections, the scaled
convolution expression~\eqref{conv_sum} is embedded in the analogous but more
favorably structured convolution expression
\begin{equation}\label{conv_sum_2}
  d_n = \sum_{m=-\sfrac{L}{2}}^{\sfrac{L}{2}-1} c_m b_{\beta m - \gamma
    n},\quad -L/2\leq n - n_0\leq L/2-1,
\end{equation}

Using the $\gamma$-fractional discrete Fourier transform
$C_p^{(\gamma)}$ (that is to say,
the fractional Fourier transform based on roots of unity parameter
$\gamma$ as in~\cite{BaileySwarztrauberSIAMRev:91}) together with the discrete Fourier transform $B_p$,
\begin{equation*}
    C_p^{(\gamma)} = \sum_{m=-\sfrac{L}{2}}^{\sfrac{L}{2}-1} c_m
        \e^{-\i \frac{2\pi \gamma mp}{L}}, \quad
    B_p = \sum_{m=-\sfrac{L}{2}}^{\sfrac{L}{2}-1} b_m \e^{-\i \frac{2\pi m p}{L}},
\end{equation*}
an application of the convolution theorem yields~\cite{Nascov:09}
\begin{equation}\label{approx_conv_sum}
    d_n = \sum_{m=-\sfrac{L}{2}}^{\sfrac{L}{2}-1} c_m b_{\beta m - \gamma n} \approx
        \frac{1}{L} \sum_{p = -\sfrac{L}{2}}^{\sfrac{L}{2} - 1} C_p^{(\gamma)} B_p
        \e^{\i \frac{2\pi\beta n p}{L}},\quad -L/2\leq n - n_0\leq L/2-1,
\end{equation}
reducing, in particular, the (approximate) evaluation of the desired
values $d_n$ in~\eqref{conv_sum} to evaluation of a discrete Fourier
transform and a $\gamma$-fractional discrete Fourier transform, both
of size $L$, followed by evaluation of the $L$-term inverse
$\beta$-fractional Fourier transform on the right-hand side
of~\cref{approx_conv_sum}. The necessary discrete Fourier transform
can of course be evaluated by means of the FFT algorithm. The
fractional Fourier transforms (FRFTs) can also be accelerated on the
basis of the FFT-based fractional Fourier transform algorithms, at an
$\mathcal{O}(L\log L)$ cost of approximately four times that of an
$L$-point FFT; see~\cite{BaileySwarztrauberSIAMRev:91}. The error
inherent in the approximation~\cref{approx_conv_sum} is a quantity of
order $\mathcal{O}(L^{-2})$, which, in our applications, generally
yields any desired accuracy in very fast computing times by selecting
appropriate values of the parameter $L$.

{\color{red}
\begin{table}
    {\footnotesize \color{red}
    \begin{minipage}[b]{0.49\hsize}\centering
        {\renewcommand{\arraystretch}{1.3}
        \begin{tabular}{|c|c|c|c|}
          \hline
          {$M$} & {Direct (s)} & {Fast (s)} & {$\varepsilon^\mathrm{Fast}$}\\\hline
          $10^1$  & $5.6\cdot 10^{-2}$  & $8.3\cdot 10^{-3}$  & $6.6\cdot 10^{-3}$ \\\hline
          $10^2$  & $8.3\cdot 10^{-2}$  & $7.9\cdot 10^{-3}$  & $1.9\cdot 10^{-7}$ \\\hline
          $10^3$  & $1.8\cdot 10^{-1}$  & $7.8\cdot 10^{-3}$  & $1.6\cdot 10^{-8}$ \\\hline
          $10^4$  & $1.4\cdot 10^{ 0}$  & $8.2\cdot 10^{-3}$  & $7.9\cdot 10^{-7}$ \\\hline
        \end{tabular}}
    \end{minipage}
    \hfill
    \begin{minipage}[b]{0.49\hsize}\centering
        {\renewcommand{\arraystretch}{1.3}
        \begin{tabular}{|c|c|c|c|}
            \hline
            {$N$} & {Direct (s)} & {Fast (s)} & {$\varepsilon^\mathrm{Fast}$}\\\hline
            $10^1$  & $1.5\cdot 10^{-3}$  & $1.1\cdot 10^{-2}$ & $6.2\cdot10^{-6}$ \\\hline
            $10^2$ & $1.2\cdot 10^{-2}$ & $7.5\cdot 10^{-3}$ & $5.8\cdot10^{-6}$ \\\hline
            $10^3$ & $8.8\cdot 10^{-2}$ & $8.1\cdot 10^{-3}$ & $4.7\cdot10^{-6}$ \\\hline
            $10^4$ & $7.1\cdot 10^{-1}$ & $8.0\cdot 10^{-3}$ & $2.6\cdot10^{-7}$ \\\hline
            $10^5$ & $7.5\cdot 10^{ 0}$ & $9.3\cdot 10^{-2}$ & $2.1\cdot10^{-9}$ \\\hline
            $10^6$ & $9.4\cdot 10^{ 1}$ & $1.5\cdot 10^{ 0}$ & $2.6\cdot10^{-10}$ \\\hline
            $10^7$ & $2.8\cdot 10^{ 3}$ & $2.3\cdot 10^{ 1}$ & $5.5\cdot10^{-9}$ \\\hline
        \end{tabular}}
    \end{minipage}
    \caption{\color{red}Computing times required for evaluation of the size-$M$ scaled convolution~\cref{discConv} by means of the Direct and Fast algorithms described in the text at a number $N$ of time points $t_n$, and errors $\varepsilon^\mathrm{Fast}$ associated with the Fast algorithm. (By definition, the Direct algorithm provides the exact convolution results, up to roundoff.) Left: $N = 10^4$. Right: $M = 5000$.}
    \label{tab:fast_conv_table}}
\end{table}
To demonstrate the accelerated scaled-convolution algorithm we
evaluate the transform~\eqref{conv_sum} for several values of $N$,
with certain coefficients $c_m$ ($-M/2\leq m \leq M/2-1$) and with
$b_q$ as in~\eqref{discConv}. (The particular selection of the
coefficients $c_m$ is immaterial in the context of the present
demonstration, but, for reference, we mention that the coefficients
used in the example were obtained as the coefficients of the $M$-term
FC expansion~\cref{fourier_approx} with $F(\omega) = \e^{-\frac{1}{4}
  (\omega - 10)^2}\e^{-\i 8\omega}$ in the interval $[8,15]$. This
specific scaled convolution arises as the method in \cref{sec:qfsx} is
applied to the evaluation of~\cref{four_trans_general} on $[0, T]$,
with $T = N\Delta t$ and $\Delta t = 0.2$.)  Letting $\widetilde d_n$
denote the approximation of $d_n$ produced by the fast
algorithm,~\cref{tab:fast_conv_table} displays the $\ell^\infty$ error
$\varepsilon^\mathrm{Fast} = \mathrm{max}_n |d_n - \widetilde d_n|$ as well as the time required
by the fast method to produce the $M$-coefficient sum at the required
$N$ evaluation points.  The computations were performed in MATLAB on
an Intel Core i7-8650U CPU\@.
}

%--------------------------------------------------------------------------------------------------
\subsection{Non-smooth $F(\omega)$: singular quadrature for 2D low
  frequency scattering}\label{sec:sing_quad}
%--------------------------------------------------------------------------------------------------
This section concerns the evaluation of the inverse transform
in~\cref{four_trans_limited} for cases in which $F$ contains an
(integrable) singularity at $\omega = 0$.  In the context of the
proposed wave equation solver, this occurs in the evaluation
of~\cref{ukW} in the $d=2$ case (where for each spatial point
$\mathbf{r}$ we have $F(\omega) = U_k^\mathrm{slow}(\mathbf{r},
\omega) \e^{\i\omega s_k}$) since, as is
known~\cite{MacCamy:65,Werner:86}, in two dimensions the solutions to
the Helmholtz equation vary as an integrable function of $\log\omega$
which vanishes at $\omega = 0$. (Special treatments are not necessary
in the $d=3$ case, where, given incident fields with smooth
$\omega$-dependence, the Helmholtz solutions vary smoothly with
$\omega$ for all real values of $\omega$~\cite{Werner:62,Kress:79}.)

To design our quadrature rule in the non-smooth case we recall the
decomposition~\cref{fourier_sing_decomp_discussion},
\begin{equation}\label{fourier_sing_decomp}
    f(t_\ell)  = \left( \int_{-W}^{-w_c} + \int_{-w_c}^{w_c} + \int_{w_c}^W \right) F(\omega)
    \e^{-\i\omega t_\ell}\,\d\omega \eqqcolon I_-(t_\ell) + I_0(t_\ell) + I_+(t_\ell),
\end{equation}
where using the notation introduced in~\cref{four_trans_general},
$I_-(t_\ell) = I_{-W}^{-w_c}[F](t_\ell)$ and $I_+ (t_\ell)=
I_{w_c}^W[F](t_\ell)$ can be treated effectively by means of the
Fourier-based quadrature method developed in
\Cref{sec:fourQuad}. Unfortunately, an application of that
approach to $I_0(t_\ell) = I_{-w_c}^{w_c}[F](t_\ell)$ would not give
rise to high-order accuracy, in view of the slow convergence of the
Fourier expansion of $F$ in the interval $[-w_c,w_c]$---that arises
from the singularity of $F$ at $\omega = 0$. We therefore develop a
special quadrature rule for evaluation of the half-interval integral
\begin{equation}\label{singular_quad_model}
    I_0^{\omega_c}\,[F](t) = \int_0^{\omega_c} F(\omega) \e^{-\i t\omega}\,\d\omega
\end{equation}
that retains the main attractive features of the integration methods
developed in the previous section: high-order quadrature at fixed cost
for evaluation at arbitrarily large times $t$.
\begin{remark}\label{fourier_points_sing}
  As in \Cref{sec:fourQuad}, the aforementioned $t$-independent errors
  can be achieved on the basis of fixed ($t$-independent) finite set,
  which will be denoted by $\mathcal{F}^\textrm{sing}$ in the present
  context, of frequency mesh points. The procedure used here, however,
  does not rely on Fourier approximation of $F$, cf.\
  \Cref{fourier_points_smooth}.
\end{remark}

In order to evaluate the Fourier integral $I_0(t_\ell)$ at fixed cost
for arbitrarily large times $t_\ell$, despite the presence of
increasingly oscillatory behavior of the transform kernel, we rely on
a certain \emph{modified} ``Filon-Clenshaw-Curtis'' high-order
quadrature approach~\cite{Dominguez:13} for non-smooth $F(\omega)$.
The \emph{classical} Filon-Clenshaw-Curtis method~\cite{Sloan:80},
which assumes a smooth function $F$, involves replacement of $F$ by
its polynomial interpolant $Q_\mathscr{N} F$ at the Clenshaw-Curtis
points followed by exact computation of certain associated ``modified
moments'' (which are given by integrals of the Chebyshev polynomials
multiplied by the oscillatory Fourier kernel). Importantly, this
classical procedure eliminates the need to interpolate the target
transform function at large numbers of frequency points as time
increases. Additionally, on account of the selection of
Clenshaw-Curtis interpolation points, the polynomial interpolants
coincide with rapidly convergent Chebyshev approximations, and,
therefore, the integration procedure converges with high-order
accuracy. The accuracy resulting from use of a Chebyshev-based
approach, which is very high for any value of $t$, actually improves
as time increases: as shown in~\cite{Dominguez:13}, the error in the
method~\cite{Sloan:80} asymptotically decreases to zero as $t
\rightarrow \infty$.

The modified Filon-Clenshaw-Curtis method~\cite{Dominguez:13} we use
in the present non-smooth-$F$ case (where $F$ is singular at $\omega
=0$ only) proceeds on the basis of a graded set
\begin{equation}\label{gradedmesh}
    \Pi_{\mathscr{M},q} \coloneqq \left\{ \mu_j \coloneqq \omega_c \left(
    \frac{j}{\mathscr{M}}\right)^q : j = 1, \ldots,
    \mathscr{M}\right\},
\end{equation}
of points in $(0,\omega_c]$ which are used to form subintervals
$(\mu_j,\mu_{j+1})$ ($1\leq j\leq \mathscr{M}-1$).  For a given
meshsize $\mathscr{N}$, each one of these subintervals is then
discretized by means of a Clenshaw-Curtis mesh containing
$\mathscr{N}$ points, and all of these meshes are combined in a single
mesh set $\mathcal{F}^\textrm{sing}$ (which contains a total of
$|\mathcal{F}^\textrm{sing}| = 2(\mathscr{M}-1)\mathscr{N}$ points)
that is to be used for evaluation of the integral $I_0$. Using this
mesh, the Clenshaw-Curtis quadrature rule is applied to the evaluation
of $I_{\mu_j}^{\mu_{j+1}}\,[F](t)$
(see~\Cref{four_trans_general}). The integral $I_0^{\omega_c}\,[F](t)$
is finally approximated by a composite quadrature rule that mirrors
the exact relation
\begin{equation}\label{gradedmeshIntRule}
    I_0^{\omega_c}\,[F](t) = \sum_{j=2}^\mathscr{M} I_{\mu_{j-1}}^{\mu_j}\,[F](t).
\end{equation}

The error introduced by this quadrature rule is discussed extensively
in \cite{Dominguez:13}, and is of course dependent on the strength of
the singularity. Briefly, in our context, and assuming $q > \mathscr{N} + 1$,
the convergence order as $\mathscr{N}\to\infty$ is determined by the number $\mathscr{M}$
of integration subintervals used: letting $I_\mathscr{N} \approx I$ denote the
approximate value produced by the composite quadrature rule using $\mathscr{N}$
Clenshaw-Curtis points per subinterval, we find the error in $I_\mathscr{N}$
satisfies~\cite[Thm.~3.6]{Dominguez:13}
\begin{equation*}
    |I[F] - I_\mathscr{N}[F]| = \mathcal{O}(\mathscr{M}^{-(1 +
    \mathscr{N})})\quad\mbox{as} \quad \mathscr{N}\to\infty.
\end{equation*}
Whenever necessary (i.e.\ for two-dimensional problems containing
nonzero content at zero-frequency), the numerical results presented in
\Cref{sec:numerical_results} were produced using the values
$\mathscr{M} = 4$, $\mathscr{N} = 8$ and $q = 9.1 > \mathscr{N}+1$.
Of course, two dimensional problems whose frequency spectrum is
bounded away from the origin, and three-dimensional problems (which
always enjoy a smooth frequency dependence even around $\omega = 0$),
do not require the use of the quadrature rule described above. The
computational cost of this algorithm does not grow with increasing
evaluation time $t$, consistent with the $\mathcal{O}(1)$ large time
sampling cost for the overall hybrid method.

%-------------------------------------------------------------------------------
%-------------------------------------------------------------------------------
\section{Fast-hybrid wave equation solver: overall algorithm
  description}\label{sec:overall_algorithm}
%-------------------------------------------------------------------------------
%-------------------------------------------------------------------------------
Utilizing a number of concepts presented in the previous sections and
additional notations, including
\begin{itemize}
\item[--] An incident field $b$ of the form~\cref{single_incid_field}
  for a given direction $\mathbf{p}$;
\item[--] A set $\mathcal{F}=\{\omega_1,\dots,\omega_J\}$ of
  frequencies (n.b.\ $\mathcal{F} =
  \mathcal{F}^\textrm{smooth}\cup\mathcal{F}^\textrm{sing}$ in the 2D
  case, and $\mathcal{F} = \mathcal{F}^\textrm{smooth}$ in the 3D
  case, cf.\ \Cref{fourier_points_smooth} and
  \Cref{fourier_points_sing})) used to discretize both the slow
  $H$-windowed Fourier transform $B_k^\textit{slow}$ (cf.\
  \cref{Bkslow_eqn}) and the corresponding slow frequency scattered
  fields $U_k^\textit{slow}$ (cf.\ \cref{Ukslow_repr});
    \item[--] A set $\mathcal{C}$ (of cardinality $N_\Gamma$) of
        scattering-boundary discretization points;
    \item[--] Sets $\mathcal{R}$ (of cardinality $N_{\mathbf{r}}$) and
        $\mathcal{T} = \left\{t_\ell:\,1 \le \ell \le N_t\right\})$ (of
        cardinality $N_t$) of discrete spatial and temporal observation points
        at which the scattered field is to be produced;
\end{itemize}
the single-incidence (see \Cref{rem:pwe})) time-domain algorithm
introduced in this contribution is summarized in the following prescriptions.
\begin{enumerate}
\item[F1] Evaluate numerically the windowed incident-field
  {\color{purple} signal} functions {\color{purple}$w_k(t)a(t)$},
  in~\cref{ak_signal_def} ($k=1,\dots, K$), over a temporal mesh
  adequate for evaluation of the Fourier transforms mentioned in
  Step~[F2].
      \item[F2] Obtain the boundary condition functions
        $B_{k}^\textit{slow}(\omega)$ at frequency mesh values
        $\omega = \omega_j\in \mathcal{F}$ ($1\leq j\leq J$) by
        Fourier transformation of the windowed signals in~[F1], in
        accordance to \cref{Bkslow_eqn}.
\item[F3] Solve a total of $J$ integral equations~\cref{CFIE_direct}
  under plane-wave incidence with incidence vector $\mathbf{p}$ (see
  \Cref{rem:pwe}) at the frequencies $\omega_j\in \mathcal{F}$, to
  produce, for each $j$, boundary integral densities $\psi^t =
  \psi^t_\mathbf{p}(\mathbf{r}', \omega_j)$, $\mathbf{r}' \in
  \mathcal{C}$.
\item[F4] For each partition index $k=1,\dots, K$ produce the frequency-domain
  scattering boundary integral density $\psi_k^\textit{slow}$ with
  support in $[-W,W]$ on the basis of the densities
  $\psi^t_\mathbf{p}$ via an application of \Cref{psikslow_repr}.
\item[F5] Complete the frequency domain portion of the algorithm by
  evaluating, at each point $\mathbf{r} \in \mathcal{R}$, the
  frequency-domain solution $U_k^\textit{slow}(\mathbf{r}, \omega_j)$
  in \Cref{Ukslow_repr} by numerical evaluation of the layer potential
  integral in that equation, using the density values
  $\psi_k^\textit{slow}(\mathbf{r}', \omega_j)$ at boundary points
  $\mathbf{r}'\in \mathcal{C}$.
\end{enumerate}
In order to evaluate the solution $u$ for all points in the set
$\mathcal{R}$, and for all times in the set $\mathcal{T}$, the
algorithm proceeds by transforming each windowed solution back to the
time domain using the quadrature methods presented in
\cref{sec:fast_four_algs}. The following prescriptions thus complete
the overall hybrid solver.
\begin{enumerate}
    \item[T0] For $k=1$ to $K$ and For each $\mathbf{r}\in\mathcal{R}$ do:
    \item[T1]
        \begin{enumerate}
        \item (3D case) Obtain
          the coefficients $c_m=c_m(\mathbf{r})$ of the Fourier series
          expansions of the form~\cref{fourier_approx} for the functions $F(\omega)
          = U_k^\textit{slow}(\mathbf{r}, \omega)$ in the interval
          $\omega\in [-W,W]$.
        \item (2D case) Obtain
          the coefficients $c_m^{(1)}=c_m^{(1)}(\mathbf{r})$ and
          $c_m^{(2)}=c_m^{(2)}(\mathbf{r})$ of the Fourier series
          expansions of the form~\cref{fourier_approx} for the
          functions $F(\omega) = U_k^\textit{slow}(\mathbf{r},
          \omega)$ in the domains $[-W,-\omega_c]$ and $[\omega_c,
          W]$, respectively. (n.b.\  $\mathcal{F}^\textrm{smooth}$
          is a discretization of the set $[-W,-\omega_c]\cup
          [\omega_c, W]$.)
    \end{enumerate}
    \item[T2]
        \begin{enumerate}
            \item (3D case) Evaluate the discrete scaled convolution
                using the fast algorithms described in
                \Cref{sec:disc_conv_fast} with coefficients $c_m = c_m(\mathbf{r})$ obtained in (T1a)
                which, on account of
                Equations~\cref{ukW} and~\cref{four_trans_general}, yields
                $u_k(\mathbf{r}, t)$ for $t \in \mathcal{T}$.
            \item (2D case)
                Evaluate two discrete scaled convolutions using the
                fast algorithms described in \Cref{sec:disc_conv_fast}
                with coefficients $c_m = c_m^{(1)}(\mathbf{r})$ and
                $c_m = c_m^{(2)}(\mathbf{r})$ to produce, for all $t\in\mathcal{T}$, $I_- = I_-(t)$ and
                $I_+ = I_+(t)$ for $F = U_k^\textit{slow}$ as in
                \Cref{sec:sing_quad}.
              \item (2D case continued) Evaluate the singular
                integral approximation $I_0$ using the methods in
                \Cref{sec:sing_quad} with the frequency points in
                $\mathcal{F}^\textrm{sing}$.
              \item (2D case continued) Evaluate $u_k(\mathbf{r}, t) =
                I_-(t) + I_0(t) + I_+(t)$ for $t \in \mathcal{T}$ (cf. \cref{fourier_sing_decomp}).
        \end{enumerate}
    \item[T4] End do
    \item[T5] Evaluate $u=\sum_{k=1}^K u_k(\mathbf{r}, t)$
        (cf.\ \Cref{uk_def}, \cref{final_u}).
    \item[T6] End
\end{enumerate}

\begin{remark}
  Calling $u_k^{W,J}$ the numerical approximations to the functions
  $u_k$ produced under the finite bandwidth $W$ and on the basis of
  the $J$ quadrature points in $\mathcal{F}$, the equation in
  algorithm step [T5] can more precisely be expressed in the form
    \begin{equation}\label{final_u}
      u(\mathbf{r}, t) \approx \sum_{k=1}^K u_k^{W,J}(\mathbf{r}, t).
    \end{equation}
    The errors $e = e(W, J)$ inherent in this approximation decay
    superalgebraically fast {\em uniformly in $\mathbf{r}$ and $t$}\/
    as $W$ grows (see \Cref{rem:1a}). {\color{orange}The
      frequency-quadrature errors resulting from the methodology
      described in \Cref{sec:fast_four_algs} for the integral in
      \cref{ukW}, further, decay superalgebraically fast (or, in the
      two-dimensional case, with prescribed high-order) as $J$
      increases, {\em uniformly in $t$}\/ (see \Cref{sec:qfsx} and
      \Cref{sec:sing_quad} for a full discussion of
      frequency-quadrature errors). As discussed in Part~II, solutions
      with quadrature errors uniform in $\mathbf{r}$ can be obtained
      either on the basis of the time-domain single layer potential
      for~\cref{w_eq} (Kirchhoff formula) for the time-dependent
      density $\psi_k$, or by means of an adequate treatment of the
      high-frequency oscillations in frequency-domain space that, in
      accordance with~\Cref{Ukslow_repr}, arise for large values of
      $|\mathbf{r}|$.}
\end{remark}

%\begin{remark}\label{rem:2}
%  As discussed extensively in Part II, a fixed (geometry-dependent)
%  number $M$ of solutions $u_k$ need to be included in the
%  sum~\cref{uk_def} to evaluate $u$ for any given time $t$ within a
%  prescribed numerical accuracy, irrespective of the time-duration for
%  which the incident signal is nonzero.  In particular, it is not
%  necessary to include an unbounded, growing number of solutions $u_k$
%  in~\cref{uk_def} as $T$ (and thus also $K$) tend to infinity. The
%  geometry-dependence of the parameter $M$ relates closely to the
%  trapping character~\cite{Morawetz:75,Morawetz:77} of the underlying
%  scattering geometry.
%\end{remark}

%-------------------------------------------------------------------------------
%-------------------------------------------------------------------------------
\section{Numerical Results}\label{sec:numerical_results}
%-------------------------------------------------------------------------------
%-------------------------------------------------------------------------------
After a brief demonstration of the proposed quadrature rule in a
simple context (\Cref{sec:four_quad_results}), this section
demonstrates the convergence of the overall algorithm
(\Cref{sec:convergence}) and it presents solutions produced by the
solver in the two- and three-dimensional contexts
(\crefrange{sec:full_solver_results}{sec:3d_scatter_results}). In
particular, \Cref{sec:full_solver_results} presents a few spatial
screenshots of long-time propagation experiments (enabled by the
time-partitioning methodology described in \Cref{sec:time_partition},
see \Cref{fig:CKK_plane_wave_offangle_tracking}), as well as, in
\Cref{fig:WhispGall_plane_wave}, results for a configuration which
gives rise to significant numbers of multiple-scattering events.
{\color{red} \Cref{sec:3d_scatter_results}, finally, presents a
  variety of three-dimensional examples, including accuracy as well as
  computational- and memory-cost comparisons with results produced by
  means of recently introduced convolution-quadrature and time-domain
  integral-equation algorithms. \Cref{sec:3d_scatter_results} also
  illustrates the applicability of the methods introduced in this
  paper to a scattering surface provided in the form of a CAD
  description (Computer Aided Design).}

%-------------------------------------------------------------------------------
\subsection{Fourier Transform Quadrature
Demonstration}\label{sec:four_quad_results}
% ------------------------------------------------------------------------------
\begin{figure}[H]
    \centering
    \begin{subfigure}[b]{0.40\textwidth}
        \includegraphics[width=\textwidth]{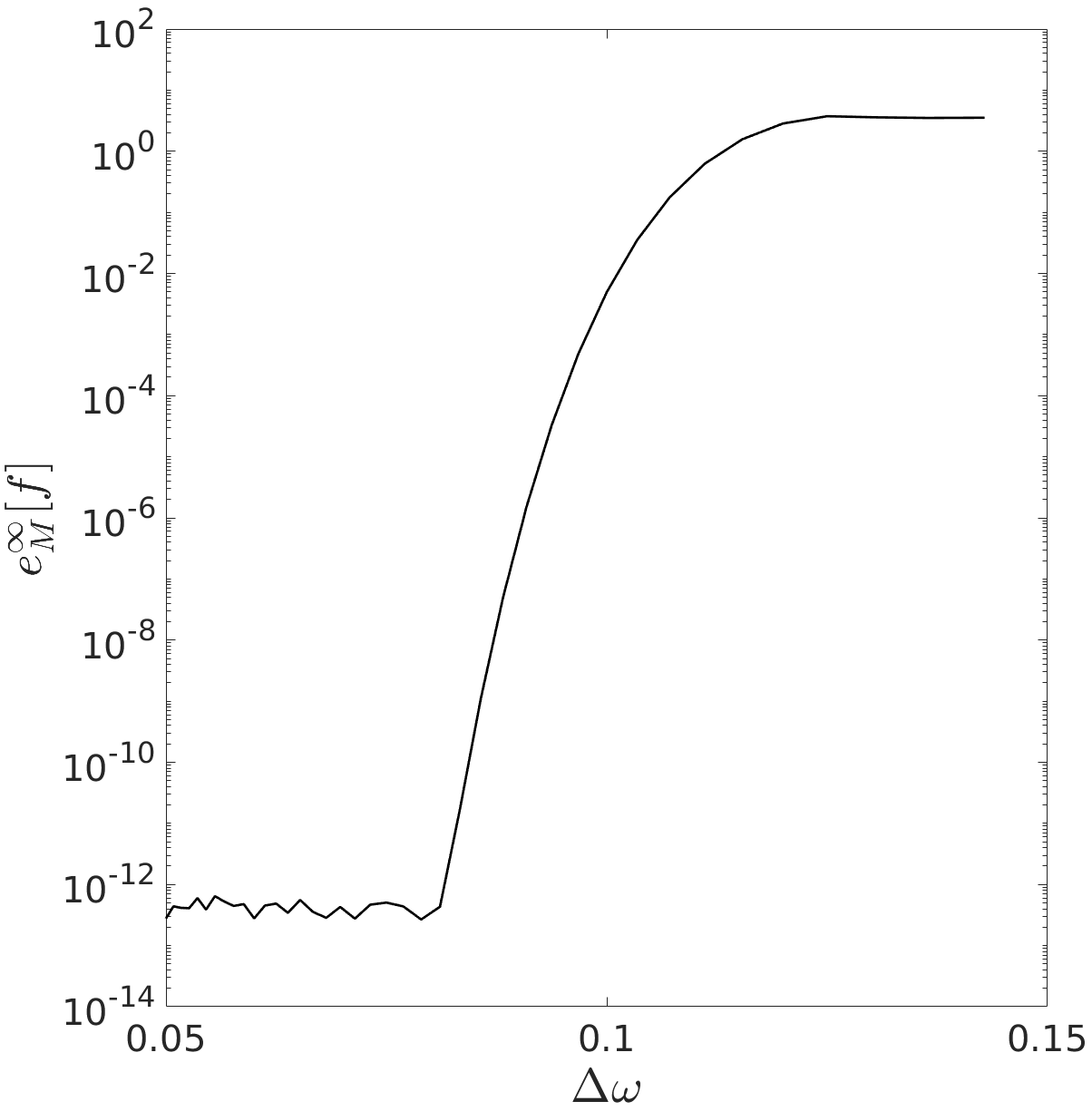}
    \end{subfigure}
    \begin{subfigure}[b]{0.40\textwidth}
        \includegraphics[width=\textwidth]{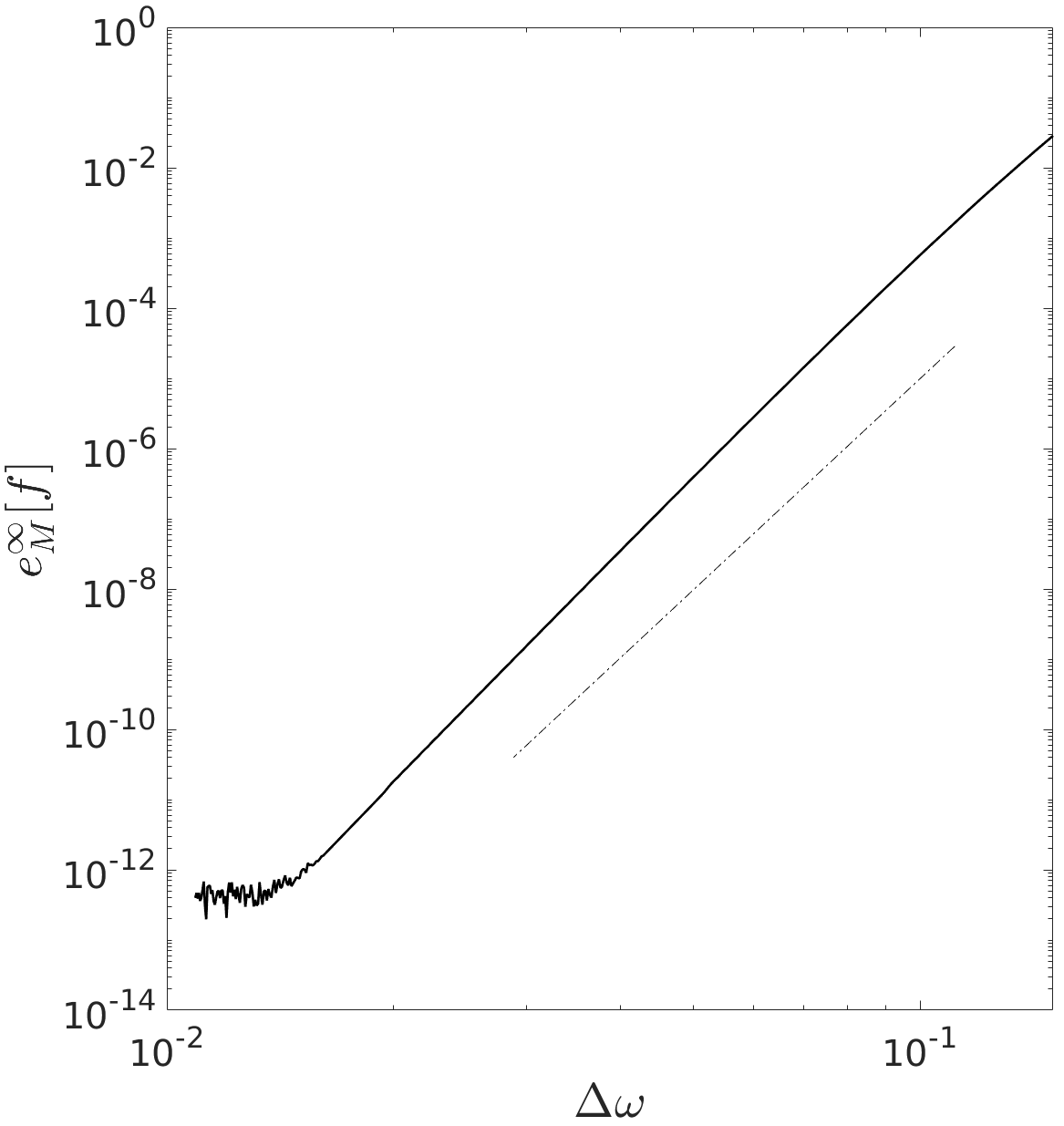}
    \end{subfigure}
    \caption{Error $e_M^\infty$ resulting from the DFT-based (left)
      and $10\textrm{th}$-order FC(Gram)-based (right)
      FRFT-accelerated Fourier Transform methods
      (cf. \Cref{Fourier-FC}) as a function of $\Delta \omega$.  The
      right figure also includes a $10\textrm{th}$-order slope, for
      reference.}
    \label{fig:fourier_quadrature}
\end{figure}
\Cref{fig:fourier_quadrature} presents results of an application (to
the function $F(\omega) = \e^{-\frac{1}{4}\omega^2} \e^{\i 10\omega}$)
of two main components of the Fourier-transform algorithms described
in \Cref{sec:fourQuad}, namely the algorithms that evaluate
trigonometric expansions~\cref{fourier_approx} by means of DFT on one
hand, and on the basis of the FC(Gram) algorithm of accuracy order
$10$, on the other (cf. \Cref{Fourier-FC}). {\color{red}(FC
  expansions of order other than $10$ can of course be used, but
  order-$10$ expansions were found perfectly satisfactory in our
  contexts.)} Noting that $|F(-12)| = |F(12)| \approx
\varepsilon_{\mathrm{mach}} $ (where $\varepsilon_{\mathrm{mach}}$
denotes machine precision), the left (resp.\ right) portion of
\Cref{fig:fourier_quadrature} displays the accuracy of the
Fourier-series based (resp.\ the FC(Gram)-based) algorithm presented
in \Cref{sec:fourQuad} for evaluation of Fourier integrals of the form
\cref{four_trans_general} in interval $[-12,12]$ (resp.\ in the
interval $[0,12]$). The fast (high-order) convergence of the
quadrature method as $\Delta \omega \to 0$ that is demonstrated in the
present simple example has a significant impact on the efficiency of
the algorithm---which requires solution of an expensive integral
equation~\cref{CFIE_direct} for each frequency discretization point
$\omega_j \in \mathcal{F^{\mathrm{smooth}}}$, $\omega_{j+1} - \omega_j
= \Delta \omega$, for $j = 2, \dots, |\mathcal{F}^{\mathrm{smooth}}| =
M$ (cf.\ \Cref{fourier_points_smooth}).

%--------------------------------------------------------------------------------------------------
\subsection{Solution convergence}\label{sec:convergence}
\begin{figure}[htp]
    \centering
    \begin{subfigure}[b]{0.32\textwidth}
        \includegraphics[height=1.7in]{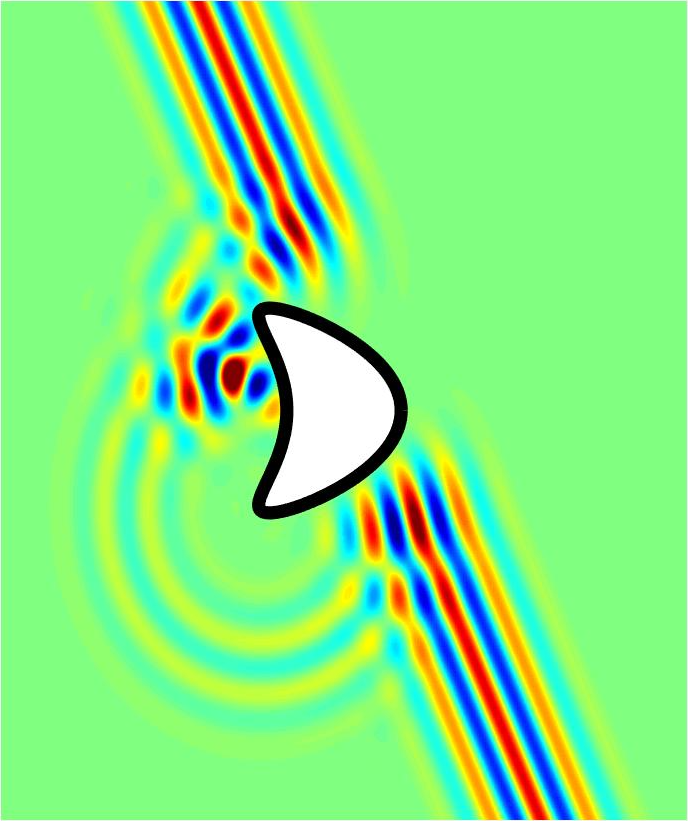}
    \end{subfigure}
    \begin{subfigure}[b]{0.32\textwidth}
        \includegraphics[width=\textwidth]{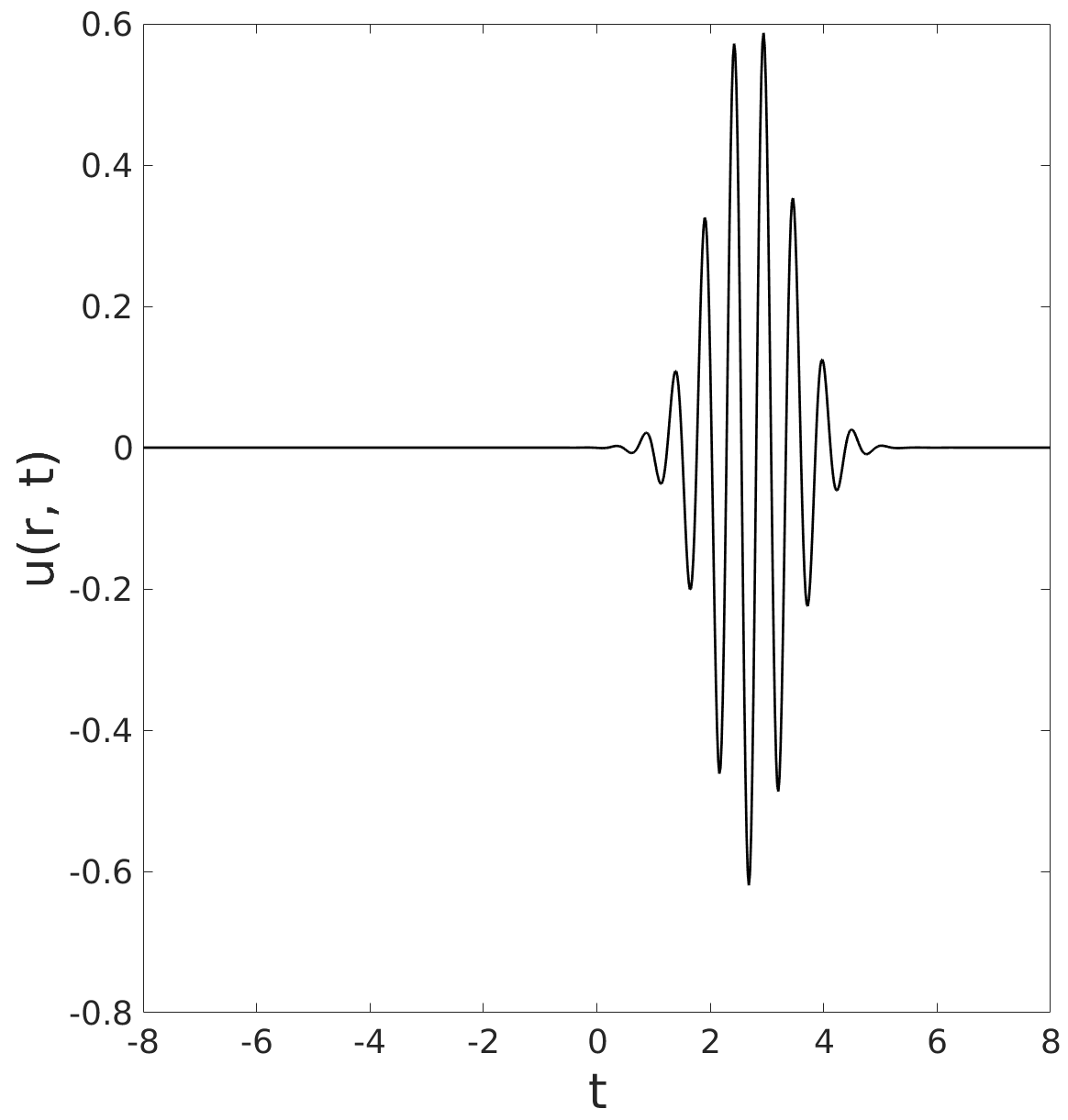}
    \end{subfigure}
    \begin{subfigure}[b]{0.32\textwidth}
        \includegraphics[width=\textwidth]{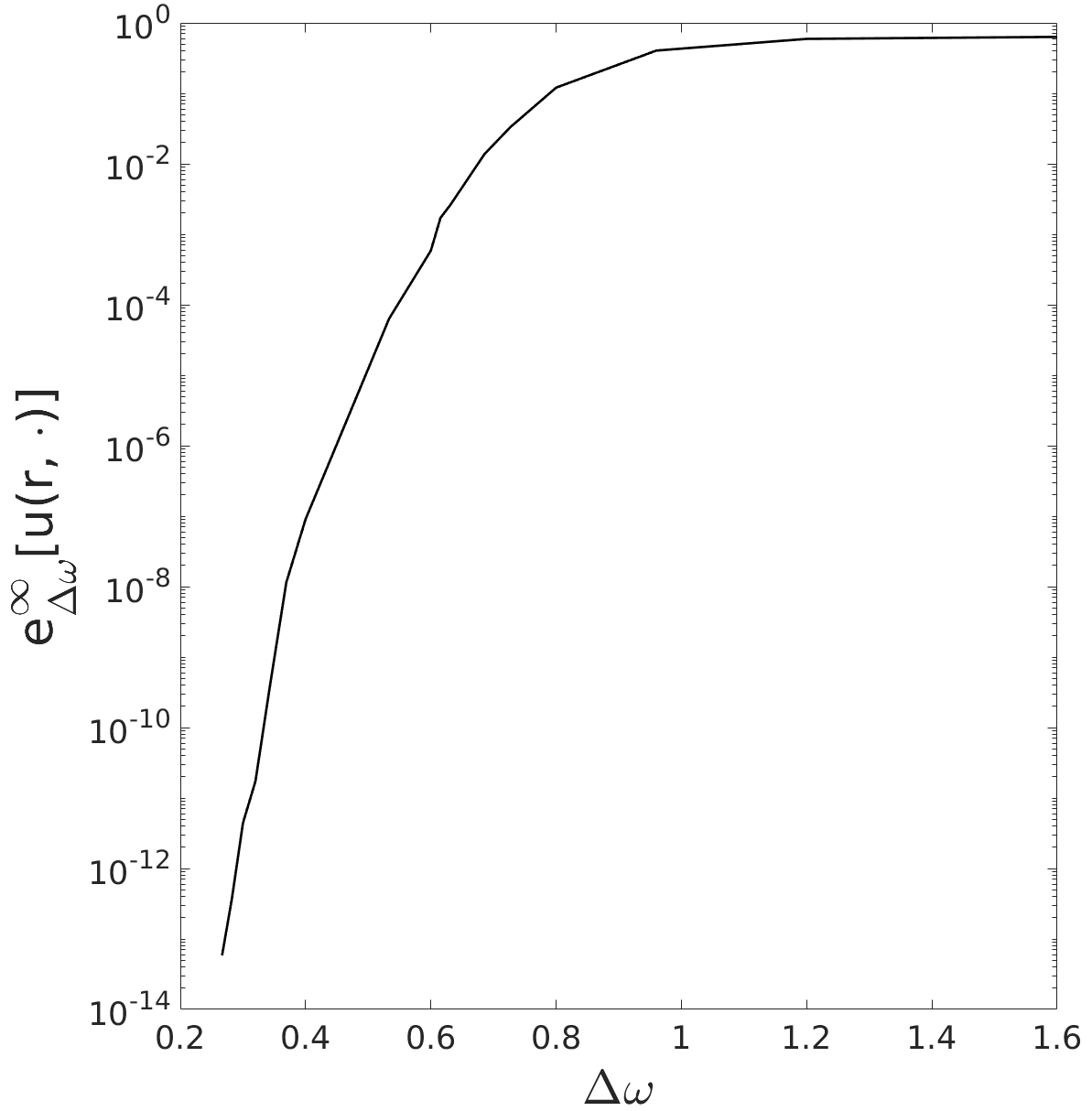}
    \end{subfigure}
    \caption{Scattered field (left), time trace  at point $\mathbf{r} = (2, 2)$
    exterior to the scatterer (center) and maximum all-time error $e_{\Delta
        \omega}^\infty$ at $\mathbf{r} = (2, 2)$ as a function of the
      frequency-domain discretization $\Delta \omega$ (right) resulting from
      an application of the overall fast hybrid method to the problem
      considered in \Cref{sec:convergence}.}
    \label{fig:conv_test}
\end{figure}

%--------------------------------------------------------------------------------------------------
This section presents solution of a problem of scattering under
incident radiation $u^\textit{inc}(\mathbf{r}, t)$ given by the
Fourier transform of the function
\begin{equation}\label{gaussian}
  U^\textit{inc}(\mathbf{r}, \omega) = \e^{-\frac{(\omega - \omega_0)^2}{\sigma^2}}
  \e^{\i\omega\mathbf{\widehat k}_\mathrm{inc}\cdot\mathbf{r}}
\end{equation}
with respect to $\omega$, with $\omega_0 = 12$, $\sigma = 2$ and,
letting $\mathbf{k} = \mathbf{e}_x + \frac{1}{2}\mathbf{e}_y$,
$\mathbf{\widehat k}_\mathrm{inc} = \frac{\mathbf{k}}{\parallel
  \mathbf{k}\parallel}$. The scatterer is a two-dimensional
kite-shaped structure $(r_1(t), r_2(t)) = (\cos(t) + 0.65\cos(2t) -
0.65,\allowbreak 1.5\sin(t))$, $(0 \leq t \leq 2\pi)$ which is also
used in the subsequent example (cf.\
\cref{fig:CKK_plane_wave_offangle_tracking}).  \Cref{fig:conv_test}
presents the time trace of the scattered field displayed in the left
image at the observation point $(2,2)$, which lies at a distance of
approximately $1.9$ spatial units from the scattering boundary. The
right image in \Cref{fig:conv_test} displays the error $e$ in the
center image as a function of $\Delta \omega$. ({\color{red} For
  simplicity, the fixed numerical bandwidth value $W = 24$ together
  with a sufficiently fine fixed spatial discretization were used in
  all cases to ensure frequency domain solution errors of the order of
  machine precision.})  The right image clearly demonstrates the
superalgebraically-fast convergence of the algorithm {\color{red}
  (relative to a converged reference solution computed with $\Delta
  \omega = 0.12$)} as the frequency-domain discretization is refined.

%--------------------------------------------------------------------------------------------------
\subsection{Full solver demonstration: 2D Examples}\label{sec:full_solver_results}
%--------------------------------------------------------------------------------------------------

This section presents results produced by the proposed methodology for
two 2D problems of sound-soft scattering---each of which demonstrates
a significant aspect of the proposed approach.

\begin{figure}[H]
\centering
\begingroup
\setlength{\tabcolsep}{3pt}
    \begin{tabular}{!{\vrule width 2pt}ScSc!{\vrule width 2pt}ScSc!{\vrule
        width 2pt}ScSc!{\vrule width 2pt}ScSc!{\vrule width 2pt}}%
        \noalign{\hrule height 2pt}%
    %\begin{tabular}{ScSc!{\vrule width 2pt}ScSc}%
    \includegraphics[width=.105\textwidth]{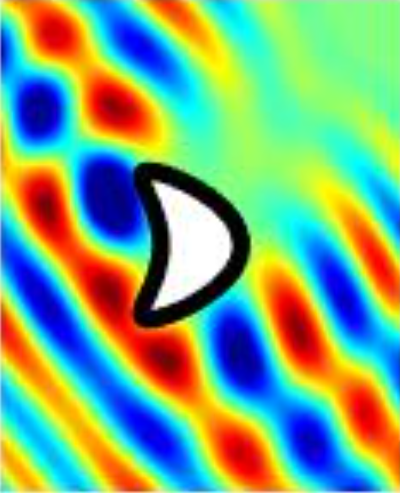}&%
    \includegraphics[width=.105\textwidth]{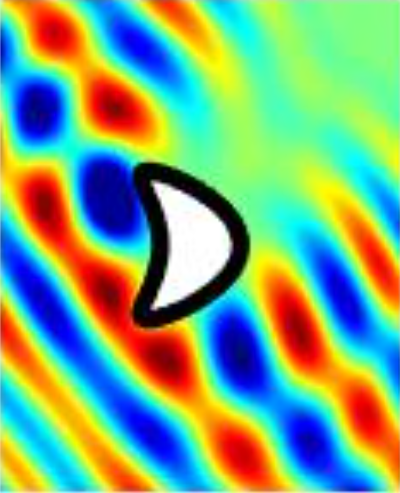}&%
    \includegraphics[width=.105\textwidth]{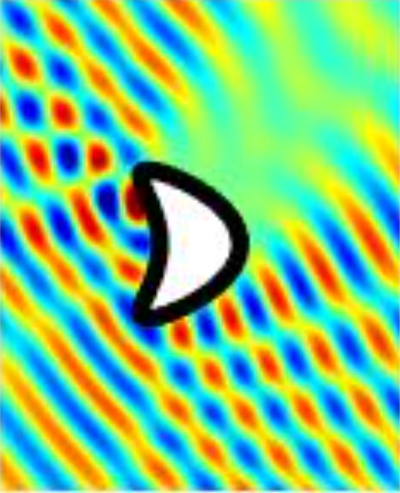}&%
    \includegraphics[width=.105\textwidth]{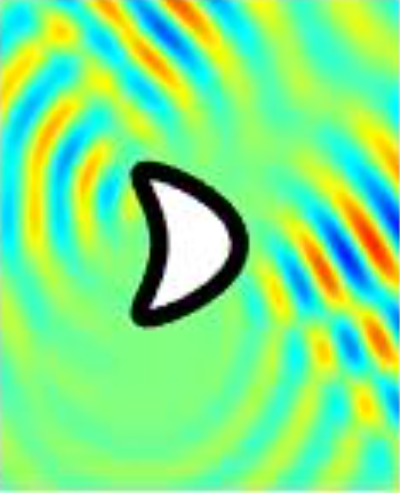}&%
    \includegraphics[width=.105\textwidth]{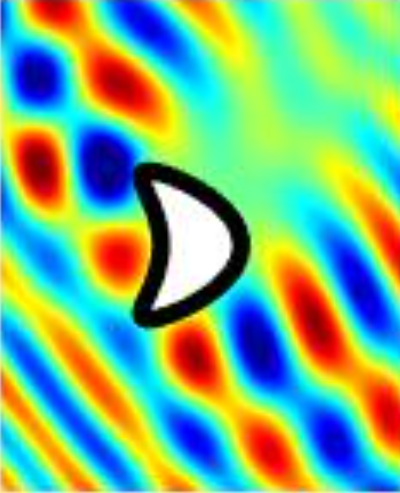}&%
    \includegraphics[width=.105\textwidth]{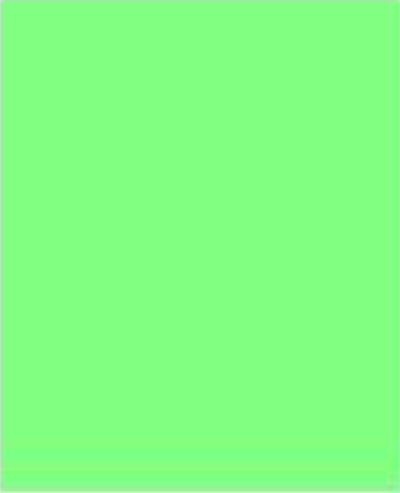}&%
    \includegraphics[width=.105\textwidth]{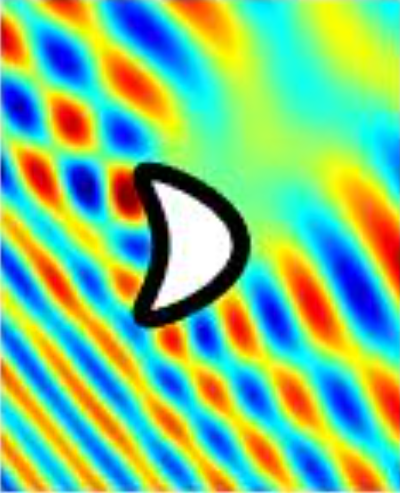}&%
    \includegraphics[width=.105\textwidth]{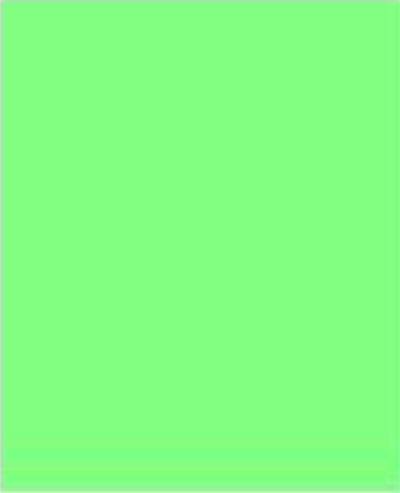}%
        \\%
    \includegraphics[width=.105\textwidth]{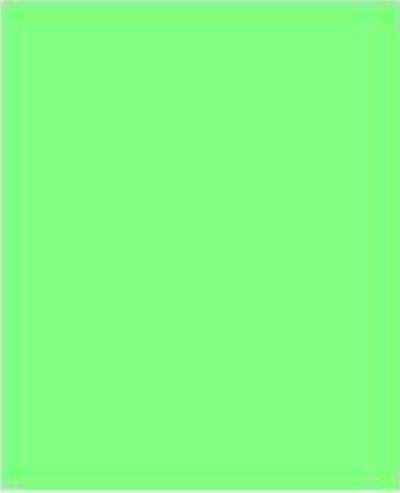}&%
    \includegraphics[width=.105\textwidth]{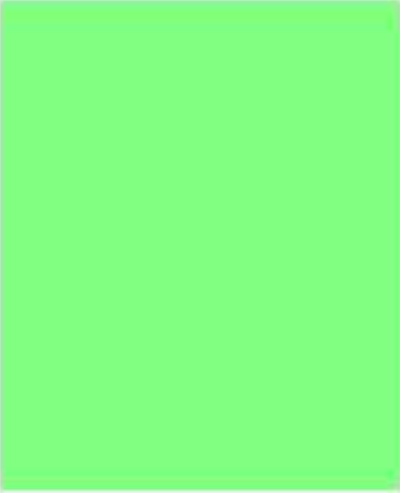}&%
    \includegraphics[width=.105\textwidth]{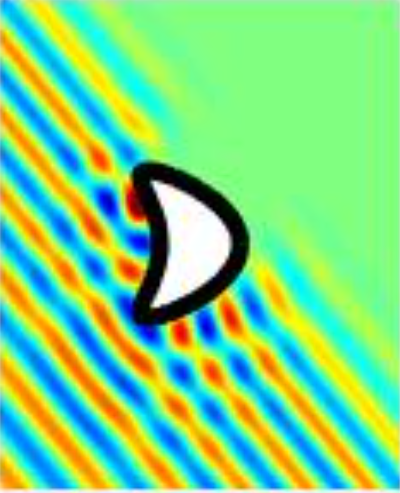}&%
    \includegraphics[width=.105\textwidth]{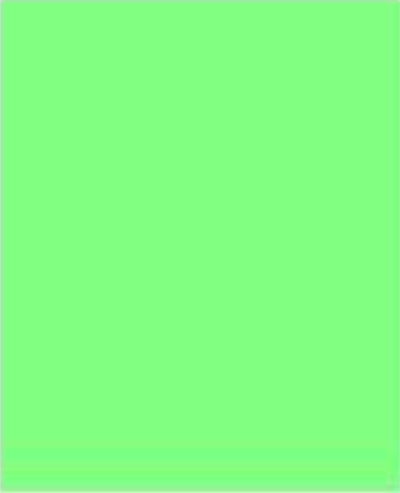}&%
    \includegraphics[width=.105\textwidth]{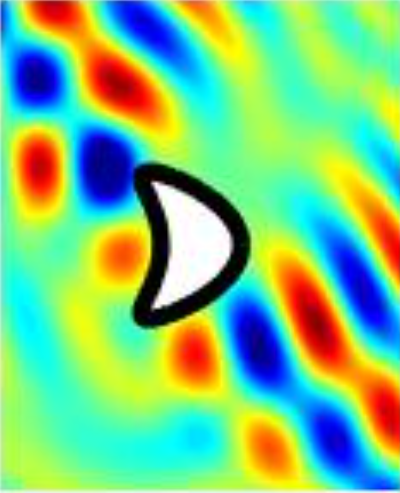}&%
    \includegraphics[width=.105\textwidth]{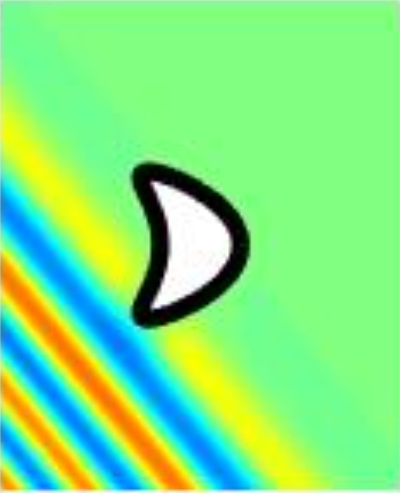}&%
    \includegraphics[width=.105\textwidth]{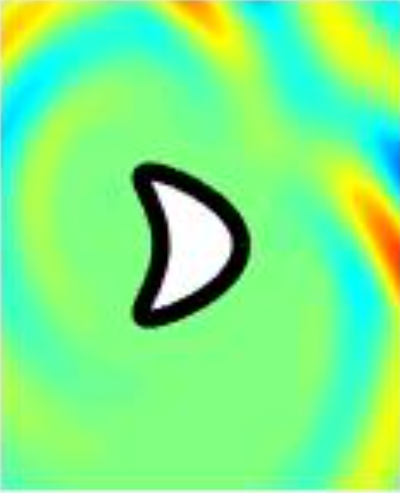}&%
    \includegraphics[width=.105\textwidth]{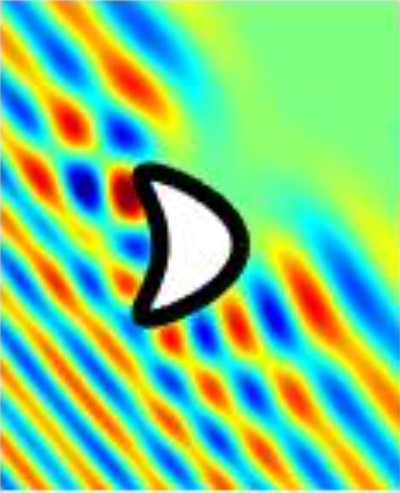}%
        \\\noalign{\hrule height 2pt}%
\end{tabular}
\endgroup
\caption{Active-partition tracking demonstration. Each of the four large panels show the solution
at increasing times, left to right.}
\label{fig:CKK_plane_wave_offangle_tracking}
\end{figure}
\begin{figure}[H]
    \centering
    \begin{subfigure}[b]{0.24\textwidth}
        \includegraphics[width=\textwidth]{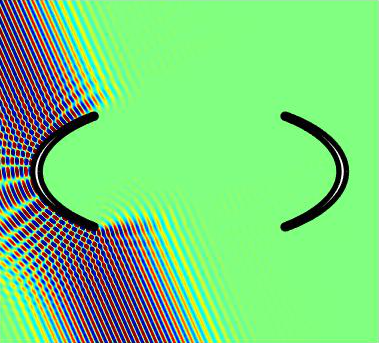}
    \end{subfigure}
    \begin{subfigure}[b]{0.24\textwidth}
        \includegraphics[width=\textwidth]{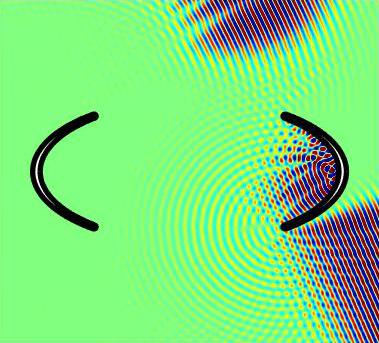}
    \end{subfigure}
    \begin{subfigure}[b]{0.24\textwidth}
        \includegraphics[width=\textwidth]{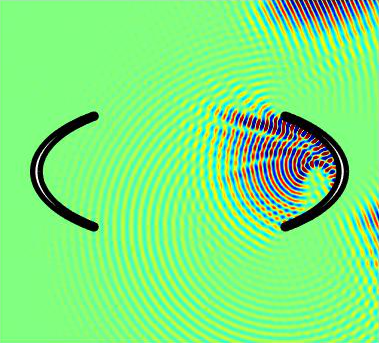}
    \end{subfigure}
    \begin{subfigure}[b]{0.24\textwidth}
        \includegraphics[width=\textwidth]{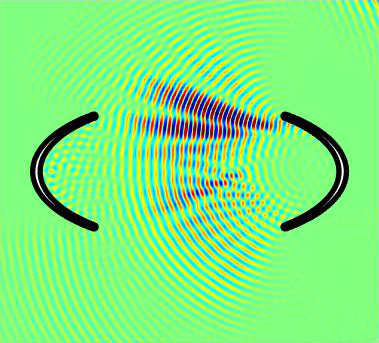}
    \end{subfigure}
    \\\vspace{.35mm}
    \centering
    \begin{subfigure}[b]{0.24\textwidth}
        \includegraphics[width=\textwidth]{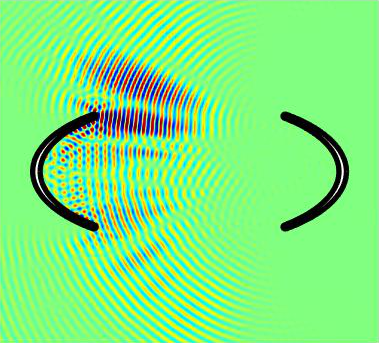}
    \end{subfigure}
    \begin{subfigure}[b]{0.24\textwidth}
        \includegraphics[width=\textwidth]{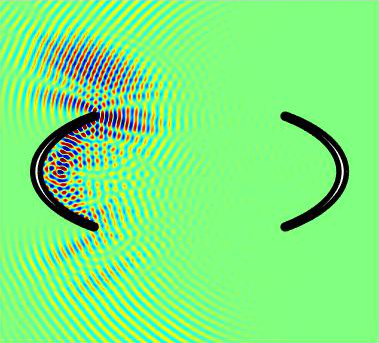}
    \end{subfigure}
    \begin{subfigure}[b]{0.24\textwidth}
        \includegraphics[width=\textwidth]{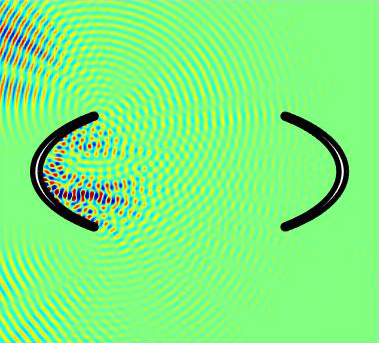}
    \end{subfigure}
    \begin{subfigure}[b]{0.24\textwidth}
        \includegraphics[width=\textwidth]{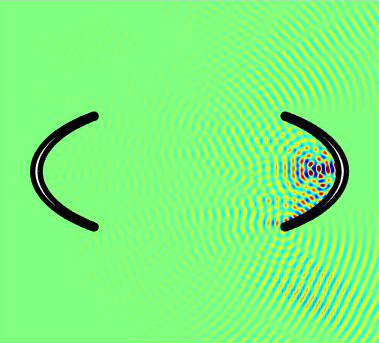}
    \end{subfigure}
    \caption{Total fields in the ``Whispering Gallery'' experiment
      mentioned in the text. Note the multiple reflections that take
      place at the elliptical surfaces which, over long propagation
      times, give rise to a significant number of scattering
      events. The time sequence starts left-to-right on the first row,
      and then continues left-to-right on the second row. }
    \label{fig:WhispGall_plane_wave}
\end{figure}
Results for incident wave-trains of longer duration which include a
time-domain chirp of the form
\begin{equation}\label{chirp_incidence}
    \begin{aligned}
        u^\textit{inc}(\mathbf{r}, t) &= {\color{purple}-a}(t - \mathbf{r} \cdot \mathbf{\widehat k}_\mathrm{inc}),\quad\mbox{with}\\
        {\color{purple}a}(t) &= \sin(g(t) + \frac{1}{4000}g^2(t)),\\
      g(t) &= 4t + 6\cos(\frac{t}{\sqrt{12}}),
    \end{aligned}
\end{equation}
and $\mathbf{\widehat k}_\mathrm{inc}$ as in \Cref{sec:convergence},
are presented in
\Cref{fig:CKK_plane_wave_offangle_tracking}, which demonstrates the time
partitioning strategy~\cref{bkdef} in conjunction with the active
partition-tracking method mentioned in \Cref{rem:tracking}. The four
large panels in this figure display results corresponding to four
subsequent time snapshots. In each one of the panels the top-left
subfigure presents the total field $u^\textit{tot}(\mathbf{r}, t)$ at
the time represented by the panel. The remaining subfigures in each
panel show the contribution to $u^\textit{tot}(\mathbf{r}, t)$ from
each one of the three corresponding time-windowed partitions used in
this example. As indicated in \Cref{rem:tracking}, blank subfigures in
\Cref{fig:CKK_plane_wave_offangle_tracking} indicate that the corresponding
partition does not contribute to $u^\textit{tot}$ at the time snapshot
represented by the panel. Using an adequate number of time windows {\color{red}
(of window width $H = 10$)} as well as a total of $200$ frequency domain
solutions {\color{red} (with bandlimit $W = 15$)}, time domain solutions at any
required time can be obtained.

\Cref{fig:WhispGall_plane_wave} demonstrates the ability of the
proposed method to account for complex multiple-scattering effects
over long periods of time. The upper left image in this figure
displays an incident wave impinging on a ``whispering gallery''
geometry; subsequent images to the right and in the lower sections of
the figure present solution snapshots at a variety of representative
times.

%--------------------------------------------------------------------------------------------------
{\color{red}
\subsection{Full solver demonstration: 3D
  Examples and Comparisons}\label{sec:3d_scatter_results}
%--------------------------------------------------------------------------------------------------

This section demonstrates the character of the proposed algorithm for
3D problems of sound-soft scattering, and it provides performance
comparisons with two solvers introduced recently.  All numerical
experiments in this section were {\color{orange} obtained by means of
  a Modern Fortran implementation of the proposed approach, using the
  Intel Fortran compiler version 17.0, on a $24$-core system
  containing two $12$-core Xeon E$5$-$2670$ CPUs.}

{\color{red}
\begin{table}
    {\footnotesize \color{red}
    \begin{minipage}[b]{0.49\hsize}\centering
        Parameter selection, comparison with~\cite{Barnett:19}.
        {\renewcommand{\arraystretch}{1.3}
        \begin{tabular}{|c|c|c|c|c|}
          \hline
            {$\omega$} & {$N$} & {$N_\beta$} & {$N_\textrm{split}$} & {$\varepsilon$}\\\hline
            $[5.5, 6.5]$  & $20$  & $150$  & $3$ & $9.0\cdot10^{-8}$\\\hline
            $[4.5, 5.5]$  & $22$  & $150$  & $2$ & $6.7\cdot10^{-9}$\\\hline
            $[3.5, 4.5]$  & $21$  & $150$  & $2$ & $1.3\cdot10^{-8}$\\\hline
            $[2.5, 3.5]$  & $21$  & $130$  & $2$ & $3.1\cdot10^{-8}$\\\hline
            $[1.5, 2.5]$  & $21$  & $120$  & $2$ & $2.7\cdot10^{-8}$\\\hline
            $[0.0, 1.5]$  & $20$  & $100$  & $2$ & $5.1\cdot10^{-8}$\\\hline
        \end{tabular}}
    \end{minipage}
    \hfill
    \begin{minipage}[b]{0.49\hsize}\centering
        Parameter selection, comparison with~\cite{Banjai:14}.
        {\renewcommand{\arraystretch}{1.3}
        \begin{tabular}{|c|c|c|c|c|}
            \hline
          {$\omega$} & {$N$} & {$N_\beta$} & {$N_\textrm{split}$} & {$\varepsilon$}\\\hline
            $[25, 45]$  & $20$  & $90$  & $3$ & $1.2\cdot10^{-4}$\\\hline
            $[20, 25]$  & $18$  & $90$  & $3$ & $1.4\cdot10^{-4}$\\\hline
            $[10, 20]$  & $15$  & $90$  & $3$ & $2.8\cdot10^{-4}$\\\hline
            $[5.0, 10]$  & $12$  & $90$  & $2$ & $1.8\cdot10^{-4}$\\\hline
            $[0.0, 5.0]$  & $10$  & $90$  & $2$ & $1.9\cdot10^{-5}$\\\hline
        \end{tabular}}
    \end{minipage}
    \caption{\color{red}
      Spatial discretizations used for the frequency-domain solver~\cite{BrunoGarza2018} in connection with comparisons with references~\cite{Barnett:19} and~\cite{Banjai:14} and associated numerical errors,
      for the $\omega$-ranges
      as listed in the first column of each table. In particular, the tables demonstrate
      that, as expected, finer discretizations need to be used, for a given
      desired accuracy, as the acoustical-size of the problems treated
      grows. In these tables, $N^2$ and $N_\textrm{split}^2$ denote the
      number of points per patch and the number of patch subdivisions of the
      original $6$-patch geometry used, respectively, so that the total number of
      degrees of freedom is $6N^2N_\textrm{split}^2$. (For the definition and
      significance of the parameter $N_\beta$, see~\cite{BrunoGarza2018}.) The
      quantity $\varepsilon$, finally, equals the numerical error at the
      spatial point $\mathbf{r} = \mathbf{r}_0$ with  $\mathbf{r}_0$ as
      indicated in the text in each comparison case, for the solution at
      frequency equal to the upper limit of the frequency interval. Mie series
      solutions were used in all cases as references for determination of the
      solution errors $\varepsilon$.
    }
    \label{tab:3D_freq_spatial_disc}}
\end{table}
}
The first example concerns a problem of scattering by a sphere of
physical radius $1.6$ (whose choice facilitates certain comparisons)
illuminated under plane-wave incidence given by
$u^\textit{inc}(\mathbf{r}, t) = {\color{purple}-a}(t - \widehat{\mathbf{k}} \cdot
\mathbf{r})$, where the signal function $T$ is given by ${\color{purple}a}(t) =
5e^{-(t-6)^2/2}$.  The frequency domain was
truncated to the interval $[-W,W]$ with numerical bandwidth $W = 6.5$,
and the problem was then discretized with respect to frequency on the
basis of $41$ ({\color{orange} $J = 80$}) equi-spaced frequencies $\omega$ in the
interval $[0, W]$. Approximate solutions to the integral equation
\cref{CFIE_direct} for each one of these frequencies were obtained by
means of the linear system solver GMRES with a relative residual
tolerance of $10^{-8}$. \Cref{tab:3D_freq_spatial_disc}
(left) lists the frequency domain spatial discretization parameters
used. The row labeled ``This work'' in
\Cref{tab:3D_barnett_comparison} presents the maximum solution error
resulting from an application of the proposed solver together with the
corresponding walltime and memory usage. We see that a computing time
of approximately four minutes and a memory allocation of 1.2 GB
suffice to produce the solution with an absolute maximum error
(measured relative to an exact solution obtained via a Mie series
representation) of the order of $10^{-7}$ at the observation point
$\mathbf{r} = (-1.8, 0, 0)$, or $0.2$ units away from the scatterer.

\begin{table}[H]
\begin{minipage}[t]{0.49\hsize}\centering
    {\footnotesize \color{red}
        {\renewcommand{\arraystretch}{1.3}
        \begin{tabular}{|c|c|c|c|}
            \hline
            {---} & {$||e||_\infty$} & {Time} & {Mem.}\\\hline
            This work & $1.6\cdot10^{-7}$  & $4.1$ & $1.2$\\\hline
            Ref.~\cite{Barnett:19} & $\approx10^{-7}$  & $101.75$  & $290$\\\hline
        \end{tabular}}
      \caption{{\color{red}Comparison with results
          in~\cite{Barnett:19}. ``This work'' data correponds to runs
          on a $24$-core computer with Sandy Bridge microarchitecture,
          while reference~\cite{Barnett:19} reports use of a $28$-core
          computer with the more recent Broadwell
          microarchitecture. The columns ``Time'' and ``Mem.'' list
          the required wall times (in minutes) and the memory usage
          (in GB).}}
    \label{tab:3D_barnett_comparison}}
\end{minipage}\hfill
\begin{minipage}[t]{0.49\hsize}\centering
{\footnotesize \color{red}
{\renewcommand{\arraystretch}{1.3}
        \begin{tabular}{|c|c|c|c|}
          \hline
            {---} & {$||e||_\infty$} & {Time} & {Mem.}\\\hline
            This work  & $2.2\cdot10^{-4}$  & $4.3$ & $1.6$\\\hline
            Ref.~\cite{Banjai:14}  & {\color{orange}$2.1\cdot10^{-3}$} & $40.1$ & $56.8$\\\hline
        \end{tabular}}
      \caption{{\color{red}Comparison with results
          in~\cite{Banjai:14}. As in that reference the computational
          times are reported in CPU core-hours.  The columns ``Time''
          and ``Mem.''  list the required CPU core-hours and the
          memory usage (in GB). The results in~\cite{Banjai:14}
          (accelerated) correspond to runs on Santa Rosa Opteron CPUs
          while the results in ``This work'' (unaccelerated) were
          obtained on the more recent Intel Sandy Bridge CPUs.}}
          %{\color{orange} (*) See discussion for determination of numerical error.} }
        \label{tab:3D_banjai_comparison}}
\end{minipage}
\end{table}

%-------------------------------------------------------------------------------

{\color{orange} This example can be related to a test case considered
in~\cite{Barnett:19}, which introduces a temporally and spatially
high-order time-domain integral equation solver, implemented in Matlab,
which relies on the built-in sparse matrix-vector multiplication
function for time-stepping and a precompiled Fortran function for
assembly of the system matrix.} A scattering configuration including a
physically realizable incident field is considered in~\cite[Sec.\
4.3]{Barnett:19} which presents computing times but which reports
errors in the median. For comparison purposes, however, it seems more
appropriate to quantify errors in some adequate norm---and, thus, we
chose to provide a comparison with results presented in Sec. 4.2 of
that paper---where a ``cruller'' scattering surface of diameter $3.2$
is used, which is illuminated by an artificial (not physically
realizable but commonly used as a test case) point source emanating
from a point interior to the surface. In lieu of solving for the
specialized cruller geometry, we compare those results to the sphere
results provided above. Noting that, with the same diameter, the
sphere has a larger surface area ($32.2$ square units) than the
cruller geometry ($23.5$ square units), the sphere problem may be
considered to be a somewhat more challenging test case in regard to
memory usage and computing cost. Errors for the test case
in~\cite[Sec.\ 4.2]{Barnett:19} can be read from the second contour
plot provided for the cruller geometry in~\cite[Fig.\ 8]{Barnett:19},
which displays an error of approximately $10^{-7}$ for $h = \Delta t
\approx 0.0367$. The memory usage and computing time required by that
test can be deduced from~\cite[Sec.\ 4.3]{Barnett:19}, and they amount
to $290$ GB of memory and $101.75$ minutes of computing time. (The
$101.75 = 23 + \frac{30}{8}\times 21$ minute computing time estimate
was obtained as the sum of precomputation and time-stepping times, as
reported in~\cite[Sec.\ 4.3]{Barnett:19}, but accounting for a
simulation over 30 time units, instead of the 8 time units reported in
that section.) In view of \Cref{tab:3D_barnett_comparison} we suggest
that, even for short propagation times, the proposed method compares
very favorably with the approach~\cite{Barnett:19} in terms of both
computational time and memory requirements.

The next example in this section concerns the scattering of a
wide-band signal of the form $u^\textit{inc}(\mathbf{r}, t) =
-0.33\sum_{i=1}^3 \exp\left[\frac{(t - \mathbf{e}_i \cdot \mathbf{r} -
    6\, \sigma - 1)^2}{\sigma^2}\right]$ from the unit sphere, where
we have set $\mathbf{e}_1 = (1, 0, 0)$, $\mathbf{e}_2 = (0, 1, 0)$,
$\mathbf{e}_3 = (0, 0, 1)$, and $\sigma = 0.1$. We solve this problem
to an absolute error level of $2.2 \cdot 10^{-4}$ (evaluated by
comparison with the exact solution obtained via a multi-incidence Mie
series representation) at the observation point $\mathbf{r} =
\mathbf{r}_0= (2.5, 0, 0)$. The frequency domain was truncated to the
interval $[-W,W]$ with numerical bandwidth $W = 45$, and the problem
was then discretized with respect to frequency on the basis of $91$
($J = 180$) equi-spaced frequencies $\omega$ in the interval $[0,
W]$. Approximate solutions to the integral equation \cref{CFIE_direct}
for each one of these frequencies were obtained by means of the linear
system solver GMRES with a relative residual tolerance of $10^{-4}$.
\Cref{tab:3D_freq_spatial_disc} (right) lists the frequency domain
spatial discretization parameters used. \Cref{fig:3D_wideband_sphere}
displays a time-trace of our solution and
\Cref{tab:3D_banjai_comparison} presents relevant performance
indicators; the core time listed for our solver was calculated as 24
times the wall time required by the parallelized frequency-domain
solver to solve all 91 frequency-domain problems in the 24-core system
used, followed by a single core run that evaluates the time
trace. Each frequency domain problem was run, in parallel, on all 24
cores.

\begin{figure}[H]
    \centering 
    \includegraphics[width=.65\textwidth]{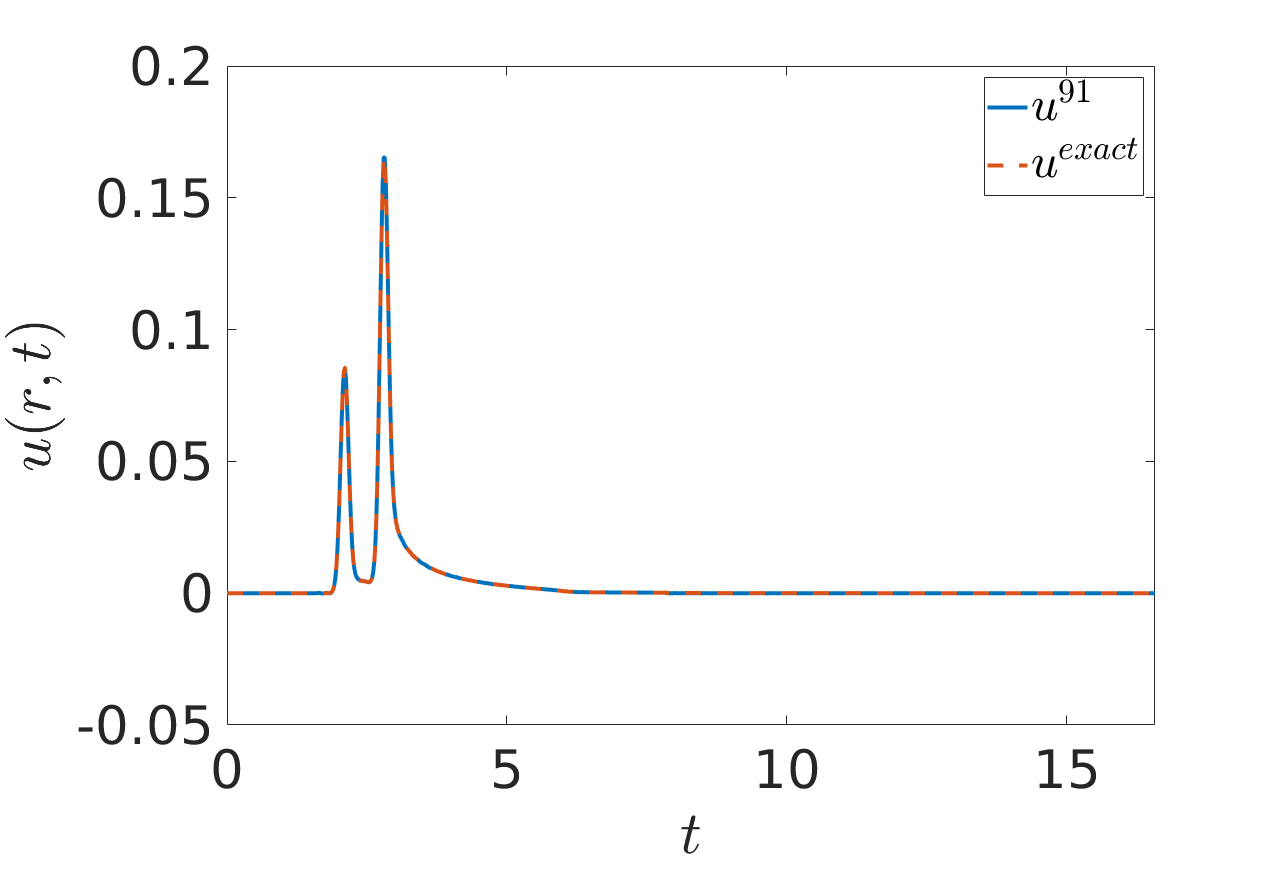}
    \caption{{\color{red}Scattering of a wide{\color{purple}-}band signal from the unit
        sphere. This figure presents the time traces of the exact and
        numerical scattered field computed (on the basis of
        {\color{orange}91}
        frequency-domain solutions) by the proposed hybrid
        method. Corresponding results produced by a novel
        convolution-quadrature implementation can be found in~\cite[Figure~5]{Banjai:14}.}}
    \label{fig:3D_wideband_sphere}
\end{figure}
This wide-band sphere-scattering example was previously considered in
the most recent high-performance implementation~\cite{Banjai:14},
including fast multipole and $\mathcal{H}$-matrix acceleration, of the
convolution quadrature method. The accuracies produced and the
computational costs required by both the present solver (for which
frequency-domain solutions were computed without use of acceleration
methods) and the one introduced in~\cite{Banjai:14} are presented in
\Cref{tab:3D_banjai_comparison}, including CPU core-hours and memory
storage. We see that, even without the significant benefits that would
arise from {\color{orange} frequency-domain} operator acceleration in
the present context, the proposed solver requires significantly
shorter computing times and lower memory allocations, by factors of
approximately ten and thirty-five, respectively, when compared with
those required by the solver~\cite{Banjai:14}. (Note: the
contribution~\cite{Banjai:14} does not directly report the solution
error, but it does provide graphical evidence by comparing the
time-traces of the solutions obtained, at the observation point
$\mathbf{r} = \mathbf{r}_0 = (2.5, 0, 0)$, by means of two different
spatio-temporal discretizations, for which it was reported that ``on
this scale [the scale of the graph] the solutions are practically
indistinguishable''. {\color{orange} A direct comparison of the
  dataset values (which was kindly provided to us by the
  authors~\cite{Kachanovska:19}) to the exact Mie solution allowed us
  to determine a maximum error of $2.1 \cdot 10^{-3}$ in the wide-band
  sphere test in~\cite{Banjai:14}.} Our solution is also graphically
indistinguishable from the time-domain response calculated by means of
a Mie series, and we report a numerical-solution error of $2.2 \cdot
10^{-4}$.)

Finally, }\Cref{fig:3d_sphere_plane_wave} presents results of an application of
the proposed algorithm to a 3D scatterer (represented by the
multi-patch CAD description displayed in the figure), for the
Gaussian-modulated incident field
\begin{equation*}
  U^\textit{inc}(\mathbf{r}, \omega) = \e^{-\frac{(\omega - \omega_0)^2}{\sigma^2}}
  \e^{\i\omega \mathbf{\widehat k}_{\mathrm{inc}} \cdot\mathbf{r}}
\end{equation*}
with $\omega_0 = 15$, $\sigma = 2$ and $\mathbf{\widehat
  k}_{\mathrm{inc}} = \mathbf{e}_z$.  {\color{red} (This figure was
  prepared using the VisIt visualization tool~\cite{VisIt}.)} A total
of $250$ frequency domain {\color{red} integral-equation solutions of
  \Cref{CFIE_direct} for frequencies below the numerical bandlimit $W
  = 25$,} which were produced by the methodology and software
described in~\cite{BrunoGarza2018}, suffice to produce the solution
for all times.
\begin{figure}[H]
    \centering
    \includegraphics[width=\textwidth]{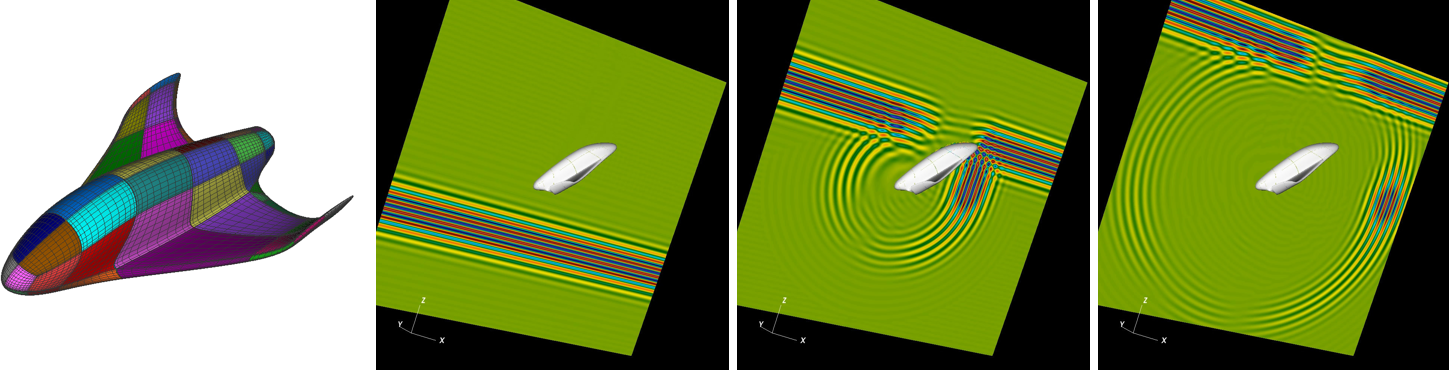}
    \caption{Field scattered by a 3D glider structure{\color{red}.}}
    \label{fig:3d_sphere_plane_wave}
\end{figure}

%--------------------------------------------------------------------------------------------------
%--------------------------------------------------------------------------------------------------
\section{Conclusions}\label{sec:conclusion}
%--------------------------------------------------------------------------------------------------
%--------------------------------------------------------------------------------------------------

This paper presents the first efficient algorithm for evaluation of
time-domain solutions, in two- and three-dimensional space, on the
basis of Fourier transformation of frequency domain solutions. The
algorithm enjoys superalgebraically fast spectral convergence in both
space and time, it runs in $\mathcal{O}(N_t)$ operations for
evaluation of the solution at $N_t$ points in time, and it can produce
arbitrarily-large time evaluation of scattered fields at
$\mathcal{O}(1)$ cost. The method is additionally embarrassingly
parallelizable in time and space, and it is amenable to
implementations involving a variety of acceleration techniques based
on high performance computing.

%--------------------------------------------------------------------------------------------------
%--------------------------------------------------------------------------------------------------
\section*{Acknowledgements}\label{sec:acknowledgements}
%--------------------------------------------------------------------------------------------------
%--------------------------------------------------------------------------------------------------
The authors gratefully acknowledge support by AFOSR, NSF and DARPA
through, respectively, contracts FA9550-15-1-0043, DMS-1714169 and
HR00111720035, and the NSSEFF Vannevar Bush Fellowship under contract
number N00014-16-1-2808. T\@. G\@. A\@. acknowledges support from the
DOE Computational Sciences Graduate Fellowship under DOE grant
DE-FG02-97ER25308. {\color{orange} Thanks are also due to Emmanuel
  Garza for facilitating the use of the existing 3D frequency-domain
  codes~\cite{BrunoGarza2018}. A number of valuable comments and
  suggestions by the reviewers are also thankfully acknowledged.}

\bibliography{tdie}

\begin{thebibliography}{10}

\bibitem{Amlani:16}
{\sc F.~Amlani and O.~P. Bruno}, {\em An {FC}-based spectral solver for
  elastodynamic problems in general three-dimensional domains}, Journal of
  Computational Physics, 307 (2016), pp.~333--354,
  \url{https://doi.org/10.1016/j.jcp.2015.11.060}.

\bibitem{Atkinson:09}
{\sc K.~E. Atkinson}, {\em The Numerical Solution of Integral Equations of the
  Second Kind (Cambridge Monographs on Applied and Computational Mathematics)},
  Cambridge University Press, 2009.

\bibitem{Babuska:97}
{\sc I.~M. Babu{\v{s}}ka and S.~A. Sauter}, {\em Is the pollution effect of the
  {FEM} avoidable for the helmholtz equation considering high wave numbers?},
  {SIAM} Journal on Numerical Analysis, 34 (1997), pp.~2392--2423,
  \url{https://doi.org/10.1137/s0036142994269186}.

\bibitem{BaileySwarztrauberSIAMRev:91}
{\sc D.~H. Bailey and P.~N. Swarztrauber}, {\em The fractional fourier
  transform and applications}, SIAM Review, 33 (1991), pp.~389--404,
  \url{https://doi.org/10.1137/1033097}.

\bibitem{BaileySwarztrauberSIAMJSciComp:94}
{\sc D.~H. Bailey and P.~N. Swarztrauber}, {\em A fast method for the numerical
  evaluation of continuous fourier and laplace transforms}, SIAM Journal on
  Scientific Computing, 15 (1994), pp.~1105--1110,
  \url{https://doi.org/10.1137/0915067}.

\bibitem{HaDuong:86}
{\sc A.~Bamberger, T.~H. Duong, and J.~C. Nedelec}, {\em Formulation
  variationnelle espace-temps pour le calcul par potentiel retard{\'{e}} de la
  diffraction d{\textquotesingle}une onde acoustique (i)}, Mathematical Methods
  in the Applied Sciences, 8 (1986), pp.~405--435,
  \url{https://doi.org/10.1002/mma.1670080127}.

\bibitem{Banjai:10}
{\sc L.~Banjai}, {\em Multistep and multistage convolution quadrature for the
  wave equation: Algorithms and experiments}, SIAM Journal on Scientific
  Computing, 32 (2010), pp.~2964--2994,
  \url{https://doi.org/10.1137/090775981}.

\bibitem{Banjai:14}
{\sc L.~Banjai and M.~Kachanovska}, {\em Fast convolution quadrature for the
  wave equation in three dimensions}, Journal of Computational Physics, 279
  (2014), pp.~103 -- 126, \url{https://doi.org/10.1016/j.jcp.2014.08.049}.

\bibitem{Banjai:11}
{\sc L.~Banjai, C.~Lubich, and J.~M. Melenk}, {\em Runge--kutta convolution
  quadrature for operators arising in wave propagation}, Numerische Mathematik,
  119 (2011), pp.~1--20, \url{https://doi.org/10.1007/s00211-011-0378-z}.

\bibitem{Banjai:09}
{\sc L.~Banjai and S.~Sauter}, {\em Rapid solution of the wave equation in
  unbounded domains}, SIAM Journal on Numerical Analysis, 47 (2009),
  pp.~227--249, \url{https://doi.org/10.1137/070690754}.

\bibitem{Banjai:12}
{\sc L.~Banjai and M.~Schanz}, {\em Wave propagation problems treated with
  convolution quadrature and bem}, in Fast Boundary Element Methods in
  Engineering and Industrial Applications, U.~Langer, M.~Schanz, O.~Steinbach,
  and W.~L. Wendland, eds., Springer Berlin Heidelberg, Berlin, Heidelberg,
  2012, pp.~145--184, \url{https://doi.org/10.1007/978-3-642-25670-7_5}.

\bibitem{Barnett:19}
{\sc A.~H. Barnett, L.~Greengard, and T.~Hagstrom}, {\em High-order
  discretization of a stable time-domain integral equation for 3d acoustic
  scattering}, 2019, \url{https://arxiv.org/abs/arXiv:1904.00076}.

\bibitem{Bayliss:80}
{\sc A.~Bayliss and E.~Turkel}, {\em Radiation boundary conditions for
  wave-like equations}, Communications on Pure and Applied Mathematics, 33
  (1980), pp.~707--725, \url{https://doi.org/10.1002/cpa.3160330603}.

\bibitem{Berenger:94}
{\sc J.-P. Berenger}, {\em A perfectly matched layer for the absorption of
  electromagnetic waves}, Journal of Computational Physics, 114 (1994),
  pp.~185--200, \url{https://doi.org/10.1006/jcph.1994.1159}.

\bibitem{Betcke:17}
{\sc T.~Betcke, N.~Salles, and W.~{\'{S}}migaj}, {\em Overresolving in the
  laplace domain for convolution quadrature methods}, {SIAM} Journal on
  Scientific Computing, 39 (2017), pp.~A188--A213,
  \url{https://doi.org/10.1137/16m106474x}.

\bibitem{Bleszynski:96}
{\sc E.~Bleszynski, M.~Bleszynski, and T.~Jaroszewicz}, {\em {AIM}: Adaptive
  integral method for solving large-scale electromagnetic scattering and
  radiation problems}, Radio Science, 31 (1996), pp.~1225--1251,
  \url{https://doi.org/10.1029/96RS02504}.

\bibitem{Borm:18}
{\sc S.~B\"{o}rm and C.~B\"{o}rst}, {\em Hybrid matrix compression for
  high-frequency problems}, 2018, \url{https://arxiv.org/abs/arXiv:1809.04384}.

\bibitem{Borm:03}
{\sc S.~B\"{o}rm, L.~Grasedyck, and W.~Hackbusch}, {\em Introduction to
  hierarchical matrices with applications}, Engineering Analysis with Boundary
  Elements, 27 (2003), pp.~405--422,
  \url{https://doi.org/10.1016/s0955-7997(02)00152-2}.

\bibitem{Bebendorf:10}
{\sc D.~Brunner, M.~Junge, P.~Rapp, M.~Bebendorf, and L.~Gaul}, {\em Comparison
  of the fast multipole method with hierarchical matrices for the
  helmholtz-bem}, Computer Modeling in Engineering \& Sciences(CMES), 58
  (2010), pp.~131--160.

\bibitem{BrunoGarza2018}
{\sc O.~P. Bruno and E.~Garza}, {\em A chebyshev-based rectangular-polar
  integral solver for scattering by general geometries described by
  non-overlapping patches}, 2018, \url{https://arxiv.org/abs/1807.01813}.

\bibitem{Bruno:01}
{\sc O.~P. Bruno and L.~A. Kunyansky}, {\em A fast, high-order algorithm for
  the solution of surface scattering problems: Basic implementation, tests, and
  applications}, Journal of Computational Physics, 169 (2001), pp.~80--110,
  \url{https://doi.org/10.1006/jcph.2001.6714}.

\bibitem{Bruno:10}
{\sc O.~P. Bruno and M.~Lyon}, {\em High-order unconditionally stable {FC}-{AD}
  solvers for general smooth domains i. basic elements}, Journal of
  Computational Physics, 229 (2010), pp.~2009--2033,
  \url{https://doi.org/10.1016/j.jcp.2009.11.020}.

\bibitem{Bruno:19}
{\sc O.~P. Bruno and A.~Pandey}, {\em Fast, higher-order direct/iterative
  hybrid solver for scattering by inhomogeneous media -- with application to
  high-frequency and discontinuous refractivity problems}, 2019,
  \url{https://arxiv.org/abs/arXiv:1907.05914}.

\bibitem{ChandlerWilde:12}
{\sc S.~N. Chandler-Wilde, I.~G. Graham, S.~Langdon, and E.~A. Spence}, {\em
  Numerical-asymptotic boundary integral methods in high-frequency acoustic
  scattering}, Acta Numerica, 21 (2012), p.~89–305,
  \url{https://doi.org/10.1017/S0962492912000037}.

\bibitem{Chen:12}
{\sc Q.~Chen, P.~Monk, X.~Wang, and D.~Weile}, {\em Analysis of convolution
  quadrature applied to the time-domain electric field integral equation},
  Communications in Computational Physics, 11 (2012), p.~383–399,
  \url{https://doi.org/10.4208/cicp.121209.111010s}.

\bibitem{VisIt}
{\sc H.~Childs, E.~Brugger, B.~Whitlock, J.~Meredith, S.~Ahern, D.~Pugmire,
  K.~Biagas, M.~Miller, C.~Harrison, G.~H. Weber, H.~Krishnan, T.~Fogal,
  A.~Sanderson, C.~Garth, E.~W. Bethel, D.~Camp, O.~R\"{u}bel, M.~Durant, J.~M.
  Favre, and P.~Navr\'{a}til}, {\em {VisIt: An End-User Tool For Visualizing
  and Analyzing Very Large Data}}, in {High Performance Visualization--Enabling
  Extreme-Scale Scientific Insight}, Oct 2012, pp.~357--372.

\bibitem{Coifman:93}
{\sc R.~Coifman, V.~Rokhlin, and S.~Wandzura}, {\em The fast multipole method
  for the wave equation: A pedestrian prescription}, IEEE Antennas and
  Propagation Magazine, 35 (1993), pp.~7--12.

\bibitem{ColtonKress}
{\sc D.~Colton and R.~Kress}, {\em Integral Equation Methods in Scattering
  Theory}, Society for Industrial and Applied Mathematics, 11 2013,
  \url{https://doi.org/10.1137/1.9781611973167}.

\bibitem{Dominguez:13}
{\sc V.~Dom\'inguez, I.~G. Graham, and T.~Kim}, {\em Filon--clenshaw--curtis
  rules for highly oscillatory integrals with algebraic singularities and
  stationary points}, SIAM Journal on Numerical Analysis, 51 (2013),
  pp.~1542--1566, \url{https://doi.org/10.1137/120884146}.

\bibitem{Douglas:93}
{\sc J.~Douglas, J.~E. Santos, D.~Sheen, and L.~S. Bennethum}, {\em {Frequency}
  {Domain} {Treatment} {Of} {One}-{Dimensional} {Scalar} {Waves}}, Mathematical
  Models and Methods in Applied Sciences, 03 (1993), pp.~171--194,
  \url{https://doi.org/10.1142/s0218202593000102}.

\bibitem{Driscoll:94}
{\sc J.~Driscoll and D.~Healy}, {\em Computing fourier transforms and
  convolutions on the 2-sphere}, Advances in Applied Mathematics, 15 (1994),
  pp.~202--250, \url{https://doi.org/10.1006/aama.1994.1008}.

\bibitem{Engquist:77}
{\sc B.~Engquist and A.~Majda}, {\em Absorbing boundary conditions for the
  numerical simulation of waves}, Mathematics of Computation, 31 (1977),
  pp.~629--629, \url{https://doi.org/10.1090/s0025-5718-1977-0436612-4}.

\bibitem{Epstein:16}
{\sc C.~L. Epstein, L.~Greengard, and T.~Hagstrom}, {\em On the stability of
  time-domain integral equations for acoustic wave propagation}, Discrete and
  Continuous Dynamical Systems, 36 (2016), pp.~4367--4382,
  \url{https://doi.org/10.3934/dcds.2016.36.4367}.

\bibitem{Frazer:84}
{\sc L.~N. Frazer and J.~F. Gettrust}, {\em On a generalization of
  {F}ilon{\textquotesingle}s method and the computation of the oscillatory
  integrals of seismology}, Geophysical Journal International, 76 (1984),
  pp.~461--481, \url{https://doi.org/10.1111/j.1365-246x.1984.tb05056.x}.

\bibitem{Gillman:15}
{\sc A.~Gillman, A.~H. Barnett, and P.-G. Martinsson}, {\em A spectrally
  accurate direct solution technique for frequency-domain scattering problems
  with variable media}, BIT Numerical Mathematics, 55 (2015), pp.~141--170.

\bibitem{Griffith:90}
{\sc J.~Griffith and M.~Nakhla}, {\em Time-domain analysis of lossy coupled
  transmission lines}, {IEEE} Transactions on Microwave Theory and Techniques,
  38 (1990), pp.~1480--1487, \url{https://doi.org/10.1109/22.58689}.

\bibitem{Grochenig:01}
{\sc K.~Gr\"{o}chenig}, {\em Foundations of Time-Frequency Analysis},
  Birkh\"{a}user Boston, 2001, \url{https://doi.org/10.1007/978-1-4612-0003-1}.

\bibitem{HaDuong:03}
{\sc T.~Ha-Duong}, {\em On retarded potential boundary integral equations and
  their discretisation}, in Topics in Computational Wave Propagation: Direct
  and Inverse Problems, M.~Ainsworth, P.~Davies, D.~Duncan, B.~Rynne, and
  P.~Martin, eds., Springer Berlin Heidelberg, Berlin, Heidelberg, 2003,
  pp.~301--336, \url{https://doi.org/10.1007/978-3-642-55483-4_8}.

\bibitem{Hagstrom:09}
{\sc T.~Hagstrom and T.~Warburton}, {\em Complete radiation boundary
  conditions: minimizing the long time error growth of local methods}, SIAM J.
  Numer. Anal., 47 (2009), pp.~3678--3704.

\bibitem{Kachanovska:19}
{\sc M.~Kachanovska}.
\newblock \emph{private communication}, ENSTA ParisTech-INRIA-CNRS, France.
\newblock 04 September 2019.

\bibitem{Kress:79}
{\sc R.~Kress}, {\em On the limiting behaviour of solutions to boundary
  integral equations associated with time harmonic wave equations for small
  frequencies}, Mathematical Methods in the Applied Sciences, 1 (1979),
  pp.~89--100, \url{https://doi.org/10.1002/mma.1670010108}.

\bibitem{Labarca:19}
{\sc I.~Labarca, L.~M. Faria, and C.~P{\'{e}}rez-Arancibia}, {\em Convolution
  quadrature methods for time-domain scattering from unbounded penetrable
  interfaces}, Proceedings of the Royal Society A: Mathematical, Physical and
  Engineering Sciences, 475 (2019), p.~20190029,
  \url{https://doi.org/10.1098/rspa.2019.0029}.

\bibitem{Lee:97}
{\sc J.-F. Lee, R.~Lee, and A.~Cangellaris}, {\em Time-domain finite-element
  methods}, {IEEE} Transactions on Antennas and Propagation, 45 (1997),
  pp.~430--442, \url{https://doi.org/10.1109/8.558658}.

\bibitem{Lin:92}
{\sc S.~Lin and E.~Kuh}, {\em Transient simulation of lossy interconnects based
  on the recursive convolution formulation}, {IEEE} Transactions on Circuits
  and Systems I: Fundamental Theory and Applications, 39 (1992), pp.~879--892,
  \url{https://doi.org/10.1109/81.199887},
  \url{https://doi.org/10.1109/81.199887}.

\bibitem{Lubich:94}
{\sc C.~Lubich}, {\em On the multistep time discretization of linear
  initial-boundary value problems and their boundary integral equations},
  Numerische Mathematik, 67 (1994), pp.~365--389,
  \url{https://doi.org/10.1007/s002110050033}.

\bibitem{MacCamy:65}
{\sc R.~C. MacCamy}, {\em Low frequency acoustic oscillations}, Quarterly of
  Applied Mathematics, 23 (1965), pp.~247--255,
  \url{http://www.jstor.org/stable/43635524}.

\bibitem{Mecocci:00}
{\sc E.~Mecocci, L.~Misici, M.~C. Recchioni, and F.~Zirilli}, {\em A new
  formalism for time-dependent wave scattering from a bounded obstacle}, The
  Journal of the Acoustical Society of America, 107 (2000), pp.~1825--1840,
  \url{https://doi.org/10.1121/1.428462}.

\bibitem{Morawetz:75}
{\sc C.~S. Morawetz}, {\em Decay for solutions of the exterior problem for the
  wave equation}, Communications on Pure and Applied Mathematics, 28 (1975),
  pp.~229--264.

\bibitem{Morawetz:77}
{\sc C.~S. Morawetz, J.~V. Ralston, and W.~A. Strauss}, {\em Decay of solutions
  of the wave equation outside nontrapping obstacles}, Communications on Pure
  and Applied Mathematics, 30 (1977), pp.~447--508.

\bibitem{Nascov:09}
{\sc V.~Nascov and P.~C. Logof\u{a}tu}, {\em Fast computation algorithm for the
  rayleigh-sommerfeld diffraction formula using a type of scaled convolution},
  Appl. Opt., 48 (2009), pp.~4310--4319,
  \url{https://doi.org/10.1364/AO.48.004310}.

\bibitem{Pedneault:17}
{\sc M.~Pedneault, C.~Turc, and Y.~Boubendir}, {\em Schur complement domain
  decomposition methods for the solution of multiple scattering problems}, IMA
  Journal of Applied Mathematics, 82 (2017), pp.~1104--1134.

\bibitem{Petropavlovsky:18}
{\sc S.~Petropavlovsky, S.~Tsynkov, and E.~Turkel}, {\em A method of boundary
  equations for unsteady hyperbolic problems in 3d}, Journal of Computational
  Physics, 365 (2018), pp.~294--323,
  \url{https://doi.org/10.1016/j.jcp.2018.03.039}.

\bibitem{Sayas:16}
{\sc F.-J. Sayas}, {\em Retarded Potentials and Time Domain Boundary Integral
  Equations}, Springer International Publishing, 2016,
  \url{https://doi.org/10.1007/978-3-319-26645-9}.

\bibitem{Sloan:80}
{\sc I.~H. Sloan and W.~E. Smith}, {\em Product integration with the
  clenshaw-curtis points: Implementation and error estimates}, Numerische
  Mathematik, 34 (1980), pp.~387--401,
  \url{https://doi.org/10.1007/bf01403676}.

\bibitem{Taflove:00}
{\sc A.~Taflove}, {\em Computational electrodynamics : the finite-difference
  time-domain method}, Artech House, Boston, 2000.

\bibitem{Werner:62}
{\sc P.~Werner}, {\em Randwertprobleme der mathematischen akustik}, Archive for
  Rational Mechanics and Analysis, 10 (1962), pp.~29--66.

\bibitem{Werner:86}
{\sc P.~Werner}, {\em Low frequency asymptotics for the reduced wave equation
  in two-dimensional exterior spaces}, Mathematical Methods in the Applied
  Sciences, 8 (1986), pp.~134--156,
  \url{https://doi.org/10.1002/mma.1670080110}.

\bibitem{Yilmaz:04}
{\sc A.~Yilmaz, J.-M. Jin, and E.~Michielssen}, {\em Time domain adaptive
  integral method for surface integral equations}, {IEEE} Transactions on
  Antennas and Propagation, 52 (2004), pp.~2692--2708,
  \url{https://doi.org/10.1109/tap.2004.834399}.

\end{thebibliography}
\appendix {\color{red} \section{General boundary {\color{orange}data}: cost
    estimates and comparisons\label{sec:pwe1}} As indicated in
  \Cref{rem:pwe}, this appendix briefly discusses possible extensions
  of the proposed hybrid method that enable solution of
  problem~\cref{w_eq} for arbitrary incidence-field functions
  $b(\mathbf{r}, t)$ on the basis of a fixed finite set of precomputed
  frequency-domain solutions that can be obtained in a reasonable
  computing time.  As suggested by that remark the approach could
  proceed via expansion of a given incident field in source- or
  scatterer-centered spherical-harmonic expansions; or
  scattering-boundary-based synthesis relying on principal-component
  analysis, etc. (For reference, note that, letting the maximum
  frequency and scatterer's physical size be denoted by $W$ and $a$,
  respectively, so that $W \cdot a$ equals the acoustical size of the
  scatterer at the maximum frequency, $\mathcal{D} =
  \mathcal{O}((Wa)^{d-1})$ ($d = 2$, $3$) denotes the number of
  frequency domain solutions required~\cite{Driscoll:94} by the
  general-incidence hybrid method at the highest frequency---which is
  also a bound on the number of frequency domain solutions required,
  per frequency, for each relevant frequency.) The resulting
  general-incidence method should prove advantageous. Discretization methods based on time-stepping may be more efficient
  than the hybrid method for small propagation times, since they do
  not require precomputations. However, the asymptotic cost estimates presented
  in what follows and the benefits arising from frequency-domain
  integral-equation acceleration, as opposed to the more complex
  time-domain integral-equation acceleration, indicate that
  significant advantages may result, even in the general incidence
  case, for scattering problems for which waves traverse the
  computational domain at least once. Of course the advantages of the
  hybrid method are much more significant in the most commonly
  considered single-incidence case---which requires a much more
  limited set of precomputed frequency-domain solutions.  An efficient
  implementation of the proposed hybrid method for general right-hand
  sides requires use of frequency-domain solvers which can rapidly
  produce multiple-incidence solutions, for each given frequency, on
  the basis of some precomputed direct matrix inverse or LU
  factorization, etc.  For sufficiently small problems a direct LU
  factorization can be used for this purpose, while, for large
  problems, fast direct
  solvers~\cite{Bebendorf:10,Borm:18,Gillman:15,Bruno:19,Pedneault:17}
  could be utilized: after a generally-significant setup cost, the
  latter methods can produce, for each relevant frequency, all
  solutions required by the proposed hybrid method at minimal
  additional cost.  To estimate the costs inherent in the use of fast
  direct methods in the context of the hybrid solver under general
  incidence, restricting attention to surface scattering problems
  considered in this paper, and letting $\mathcal{\widetilde O}(X) =
  \mathcal{O}(X \log(X))$, the $\mathcal H$-matrix setup cost for each
  one of the $\mathcal{O}(Wa)$ required frequencies is a quantity of
  the order of $\mathcal{\widetilde O}((Wa)^{d-1})$. Consideration of
  numerical experiments presented in~\cite[Tables 7 \& 8]{Borm:03}
  and~\cite[Table 6]{Bebendorf:10} suggest that, if computed by means
  of a $\mathcal{H}$-matrix approach, the cost to obtain all
  $\mathcal{D}$ frequency-domain solutions required by the hybrid
  method, for a given frequency, is itself a quantity of order
  $\mathcal{\widetilde O}((Wa)^{d-1})$.  Thus, the cost required for
  the evaluation of all necessary frequency domain solutions for all
  $\mathcal{O}(Wa)$ frequencies may be estimated as a quantity of the
  order of $\mathcal{\widetilde O}((Wa)^{d})$.  Volumetric time
  stepping methods (based, say, on finite-differences or
  finite-elements) over a domain of acoustical size $Wa$, on the other
  hand, require at least a spatio-temporal discretization, and thus a
  computing cost, of order $\mathcal{O}((Wa)^{d+1})$ (for simulations
  long enough that a single crest can traverse the complete
  computational domain) if a fixed number of points per wavelength are
  used. Even larger discretizations would be required, in addition, if
  dispersion is compensated for by decreasing the time-step $\Delta t$
  faster than the frequency grows. In any case, a clearly higher
  operation count results for volumetric solvers than the one
  presented above for the hybrid method with $\mathcal{H}$-matrix
  frequency-domain precomputation.  If appropriately-accelerated
  time-domain integral equations are used
  instead~\cite{Barnett:19,HaDuong:03,HaDuong:86}, finally, a cost of
  $\mathcal{\widetilde O}((Wa)^d)$ would result that is asymptotically
  comparable to the fast-hybrid/direct-solver approach described
  above---and the relative advantages would depend on other specific
  characteristics of the methods used, including, in particular,
  accuracy order and dispersion in time, as well as the quality of the
  respective acceleration methods used. A comparison of the
  unaccelerated hybrid method to a recent unaccelerated time-domain
  integral equation solver is presented
  in \Cref{sec:3d_scatter_results}. For most problems arising in
  applications, however, the incident fields can be represented by a small
  number of sources, and, in such cases an accelerated version of the
  hybrid method proposed in this paper already provides clear
  significant advantages, including improved cost asymptotics, over
   other methodologies considered in this section.
 }
\end{document}